\documentclass{article}
\usepackage{graphics,amssymb,amsfonts,eucal,fullpage,epsfig,graphpap,amsbsy}
\newcommand{\bA}{\mathbf{A}}
\newcommand{\bF}{\mathbf{F}}
\newcommand{\bH}{\mathbf{H}}
\newcommand{\bM}{\mathbf{M}}
\newcommand{\blf}{\mathbf{f}}
\newcommand{\bmu}{\boldsymbol{\mu}}
\newcommand{\bomega}{\boldsymbol{\omega}}
\newcommand{\bxi}{\boldsymbol{\xi}}
\newcommand{\Rset}{\mathbb{R}}
\newcommand{\bx}{{\mathbf{x}}}
\newcommand{\half}{\frac{1}{2}}
\newcommand{\Exp}{\mbox{\rm Exp}}
\newcommand{\skw}[1]{\mathrm{skew} \left[ #1\right]}
\newcommand{\cay}[1]{\mathrm{cay} \left( #1 \right)}
  \title{Geometric integration on spheres and some
  interesting applications}
  
  \author{D. Lewis $^*$ and N. Nigam $^\dag$\footnote{$^*$ Department of Mathematics, University of California, Santa Cruz, $^\dag$Department of Mathematics and Statistics, McGill
  University, Montreal }}
  \date{April 2, 2002}
  
\begin{document}
\maketitle
\begin{abstract}
Geometric integration theory can be employed
when numerically solving ODEs or PDEs with constraints.
In this paper, we present several one-step algorithms of various orders
for ODEs on a collection of spheres. To demonstrate the
versatility of these algorithms, we present representative calculations 
for reduced free rigid body motion (a conservative ODE) and a discretization of
micromagnetics (a dissipative PDE). We emphasize the role of isotropy
in geometric integration and link numerical integration schemes 
to modern differential geometry through the use of partial connection
forms; this theoretical framework generalizes moving frames and connections
on principal bundles to manifolds with nonfree actions.
\end{abstract}

  
 
 \section{Introduction}
In this article, we describe a set of algorithms on
multiple copies of $S^2$ (with possibly nonlinear interactions),
  with applications in material science. The physical processes we are
  interested in are modeled by PDEs of the form
  \begin{equation}\label{generalsystem1} \frac{\partial}{\partial
  t}\bmu({\mathbf{x}},t) = {\bf A(\bmu({\mathbf{x}}}, t))\,\times\,\bmu({\mathbf{x}},t),\qquad \bmu({\mathbf{x}},0) = \bmu_0({\mathbf{x}}).\end{equation}
  Here $\bmu: {\mathcal B}  \to \Rset^3$ is the field describing the physical process of interest, with 
  ${\mathbf{x}}\in {\mathcal B} $ denoting the spatial variable on some closed, compact subset of $\Rset^3$, and 
  ${\bf  A(\bmu)}$ is a (typically nonlinear) function of $\bmu$. 
  One immediately notices that $\|\bmu({\mathbf{x}},\cdot)\|$ is constant for
  all time; this constraint should be respected by the numerical
  methods used to study such systems. Classical numerical integrators fail to 
  preserve this norm, and this motivates our use of geometric
  integration.
  
The paper is organized as follows: we begin by introducing a problem
  in micromagnetics. We
  then describe the Lie group framework suitable for
  the general problem (\ref{generalsystem1}). The new, and fairly
  geometric, mathematical constructs of partial moving frames and partial
  connections, which can be
  used to select generators for use in numerical methods, are introduced
  in section 3. Readers who are unfamiliar with modern differential geometry, 
  particularly principal bundles, are advised to skip this
  section at a first reading. In section
  \ref{sectionalgo}, we present some one-step algorithms of different
  orders. An
  arbitrary function appears in these algorithms; different choices of this function
yield distinct discrete trajectories. 
We describe a possible choice that is related both to moving frames
 (\cite{mc1,mc2,lo})
and, in the case of a geometric version of the forward Euler method,
to discretization error minimization.
  Work is in progress to identify analogous function choices
  for higher order methods (\cite{lno}). We present the results of 
numerical experiments carried out using some of these 
functions for the symplectically reduced free rigid body, a Hamiltonian 
system on the two-sphere. The conservative nature of the rigid body system 
facilitates both the derivation of appropriate functions and the assessment 
of the relative performance of the algorithms. 
Finally, we present some representative numerical results 
for micromagnetics that indicate that these geometric integration schemes  
are competitive with conventional numerical algorithms. Our numerical
experiments were carried out using first, second, and fourth order methods.
  
 \subsection{ Why use geometric integrators? }
 
  One of the intrinsic features of the system (\ref{generalsystem1}) is
  that the vector field $\bmu(\bx,\cdot)$ evolves on a sphere. However, a classical integrator will
  update  $\bmu(\bx,t)$ at time $t$ using the
  approximation
  \[ \bmu(\bx,t+\Delta t)\approx \bmu(\bx,t) + \bF(\bmu(\bx,t), t, \Delta t).\] The particular
  form of $\bF(\cdot,\cdot,\cdot)$ depends on the algorithm chosen;
  however, it is clear that such updates correspond to {\it
  translations} of $\bmu(\bx,t)$, not {\em rotations}. Thus, a classical
  integrator does not account for the fact that $\bmu$ evolves on a sphere.
   This constraint is
  difficult to efficiently impose in practice.
  
  A naive approach is to keep track
  of changes in the norm $\|\bmu(\bx,t)\|$ during a numerical experiment, and
  renormalize the iterates after a prescribed tolerance has been
  exceeded. However, this renormalization is
  equivalent to the aphysical addition (or subtraction) of energy to
  the system and is therefore an undesirable solution. In addition,
  this renormalization would also affect ${\bf A(\bmu)}$. We shall see
  that in the context of micromagnetics, this change is nontrivial
  and nonlinear.
  
  Another observation is that the component of ${\bf A}(\bmu)$ parallel
  to $\bmu$ does not influence the solution curve $\bmu(\bx,t).$ However,
  this component does alter the discrete trajectories generated by
  numerical algorithms.  An appropriate selection of the normal component can
  improve the performance of the scheme. Preliminary numerical 
  investigations suggest that such corrections may allow increased 
  accuracy in the capture of key features, e.g. orbits and/or conserved
  quantities, at a moderate computational cost, without the introduction of
  significant numerical artifacts, e.g. numerical accelerations.
  
  \subsection{Who needs geometric integrators?}
  The philosophy that numerical algorithms should respect properties
  intrinsic to the system --- momentum conservation, evolution on a
  manifold, Hamiltonian structure --- is not a novel idea. (See, e.g.,
\cite{deVog}, \cite{dragtfinn}, \cite{veselov}). Many classes of algorithms
 have been  developed with this guiding principle,  especially for
   energy-conserving or  symplectic
 systems. As far as we know, the application of such methods in the context of numerical micromagnetics is relatively new. Encouraged by the success of these algorithms in an industrial context (\cite{jointreport}), we believe that these techniques 
 will find wider use in the material science community.
  
  Modern magnetic materials are used in an increasingly large number
  of applications, including thin film read heads and recording
  media \cite{parkin}, nanocrystalline
  permanent magnets, and magnetohydrodynamic fluids. In addition, there has been
  much interest in the use of ``smart
  materials'', including magnetostrictive actuators  and  organic
  ferromagnets \cite{iraqi}. These magnetic materials exhibit different responses
  corresponding to varying magnetic fields. For example, the resistance of a
  read-head ferromagnetic sensor changes as the device rotates over a
  recording medium. The magnetization of the material directly
  interacts with the other physical and chemical characteristics of the material; 
  micromagnetics theory describes this interaction at a microscopic level. Of
  particularly great interest is the correlation of physical microstructure
  and magnetization; the ability to predict the response of one to
  variations in the other is crucial to the further development of
  these materials and devices.

Another physical system whose mathematical model strongly resembles that of micromagnetics arises in the study of liquid crystals.
  Nematic and smectic liquid crystals form the basis of the operation of many every-day devices, such as LCD's,
  telecommunication devices, thermometers, projection systems --- even
  mood rings. The devices operate on the principle that a suitable applied
  field will change the orientation of the liquid crystals, while
  conserving the  pointwise-norm. Therefore, the resulting mathematical model has the same
  constraints as those in the micromagnetics situation. This system has been studied extensively, yet a norm-conserving algorithm has only been presented for an extremely 
simple situation (\cite{crystal}).

  The study of long-chain molecules such as those occurring in bio-molecular
systems yield another possible application area for the methods we  develop. These chemical systems have more complicated constraints on the geometry of possible configurations; the configuration sought is one that minimizes a free-energy functional. Statistical thermodynamics considerations permit the reformulation of these optimization problems as evolutions in time on given manifolds to a steady state (\cite{kollman},\cite{wang}). Such studies are used in the development of new drugs.


 \subsection{A mathematical model of micromagnetics}

  At any point in time, the magnetization $\bmu$ is  constant in  small regions in the material
  (termed a domain), and  the switching of these domains from one state to
  another is the basis of the functioning of devices built with these
  magnetic substances \cite{white}. Ideally the entire device would be one large domain, which
  would switch instantaneously when an applied field was imposed. In
  practice, there are several domains, and  the net magnetization
  in a desired direction is not optimal. We are interested in tracking
  the evolution of these domains, which entails following the local
  behavior of $\bmu$. In industrial applications, a ferromagnetic device  is
  subjected to changing magnetic fields (corresponding to the various
  applications --- a disk moving underneath a read-head, changes in applied
  voltages for magnetorheological fluids, etc). One is interested in the response of the
  device to the gradual changing of these external fields. 
  
   A model that is widely used in the industry is the {\it
  Landau-Lifshitz-Gilbert (LLG)} model of micromagnetics, 
 which describes the evolution of the state of magnetization $\bmu$ in
  a ferromagnetic sensor, occupying a region ${\mathcal B}$ in
  space. The LLG equation for the magnetization $\bmu(x,t)$  is given by 
\begin{equation} \label{LLG} \frac{\partial}{\partial t}\bmu=-\bmu\times \bH_{\rm eff}(\bmu)-\lambda \bmu\times(\bmu\times \bH_{\rm eff}(\bmu)), \qquad \|\bmu({\mathbf{x}})\|=1\,\,\,\forall \ {\mathbf{x}}\in \mathcal{B}.\end{equation}
Here $\lambda$ is a damping parameter and $\bH_{\rm eff}$ is an effective
  magnetic field, described in detail below. 
  The first term on the right hand side of (\ref{LLG}) describes the
  (undamped) {\it Larmor precession} of 
$\bmu$ about $\bH_{\rm eff}$ and is derived from first principles
  \cite{brown,aharoni}. It is observed
  in physical experiments, however, that changes in magnetization
  decay in finite time. The second term in (\ref{LLG}) is a
  phenomenological term (called the {\it Gilbert damping term}, 
see \cite{gilbert,landau}), added to describe this damping behavior; it
  cannot be derived from first principles. There are
  situations under which this system is stiff (see, for example, \cite{stiffdella,tsiantos}), and issues of numerical
  stability of the integration scheme
  are therefore important. We also thank one of the referees for pointing out \cite{donahue00}, where it is suggested that poor representations of the exchange energy may lead to the observed stiff behavior.
  
  The effective field, which causes the magnetization to change,  is derived from energy considerations
\cite{brown} and varies nonlinearly with $\bmu$. More precisely, 
 \begin{equation}\label{heff}\bH_{\rm eff}(\bmu)=A\Delta \bmu+\bmu_0
  \left(-\nabla \phi+ \bH_{\rm app} \right)+K(\bmu\cdot \mathbf{e})\mathbf{e}.
  \end{equation} 
The parameters $A$ and $K$ are material constants of the
  permalloy being studied, and $\bmu_0$ is the  permeability of
  free-space. The field $A\Delta \bmu$ is called the {\it exchange} field, 
preventing rapid spatial  variations of $\bmu$ and the formation
  of arbitrarily fine domains (see \cite{donahue97,scholz} for examples
 on how this term is computed in general). The final contribution in (\ref{heff}) is due to the
  nature of  ferromagnetic
  crystals, which causes the magnetic moments to align in preferred
  directions.  This effect is incorporated in the LLG model
  through the {\it uniaxial anisotropy} field $K(\bmu\cdot
  \mathbf{e})\mathbf{e}$.   The external
  {\it applied field} is denoted by $\bH_{\rm app}$.   The nonlinear, nonlocal contributions of
  $\bmu$ arise through the {\it demagnetizing field}, $-\nabla \phi$,
  where $\phi$  solves the Poisson
  problem with suitable boundary  \begin{equation}\label{poisson}\Delta \phi = \nabla\cdot \bmu\quad \mbox{in}\,\,\Rset^3;\,\, \left[\phi\right]=0,\quad  \left[ \frac{\partial \phi}{\partial n}\right]
  = \bmu\cdot {\mathbf n}  \,\,
  \mbox{on} \,\, \partial{\mathcal{B}},\end{equation}  and radiation conditions\[ \phi=o\left (\frac{1}{|\bx|}\right )\,\, \mbox{as} \ |\bx|\rightarrow
  \infty. \]
  Here $[u]$ denotes the jump of the function $u$ across $\partial{\mathcal B}$. 
Many different methods exist for the calculation of this field,
  including the use of the full Maxwell system, FFT techniques, finite element methods, multigrid
  approaches,  finite
  differences, and recently, fast multipole
  methods,
  \cite{fredkin1998d,oti1993,polstyanko1999,tsukerman1993,monk1999,aharoni1991,aharoni1999,asselin1986,fredkin1990,luskin1993,fidler1996,fidler1999}, among others. While there are still several unsettled issues
in the area of demagnetizing field calculations (see \cite{aharoni1999} for a sharp critique on existing
methods), it is not our intention in this project to duplicate this work. Instead, we shall focus on
developing a time-stepping method that is robust, accurate, and  requires relatively few
 field evaluations.
  
    Theoretical developments in micromagnetics are driven by 
  industrial demands, and the need for accurate  algorithms is now
  imperative. Conventional algorithms are still being employed for
   highly sensitive calculations on large sensors, and are becoming increasingly
  inadequate.  Moreover, the  time-scales inherent in these problems vary from
  nano-seconds (in disk drives) to  tens of  seconds; hence integration 
  techniques that remain effective over
  long times are required. An ideal integrator would resolve solutions
  accurately over very small time steps, while allowing large
  time steps to be taken when the system evolves over a long period.

  \subsection{Potential problems with renormalising}
  We have already described some of the issues with classical
  integrators with respect to norm conservation. Here we show that renormalizing $\bmu$ to conserve the norm while using conventional
  integrators changes the potential $\phi$ in a nonlinear fashion; 
  it is easy to construct a simple example in which the renormalization
  introduces a significant change in the demagnetizing field. 
  Assume that the magnetization $\bmu$ satisfies
  \begin{equation}\label{renorm_exam}
  \bmu(x, y, z)=a \, {\bf i} + b(x) \, {\bf j}, 
  \end{equation}
  for some constant $a$ and scalar function $b$; here 
  ${\bf i,j}$ denote the usual unit vectors in the $x$ and $y$ directions.
  Before renormalization, the field $\bM$ is divergence--free: 
  $\nabla \cdot \bmu \equiv 0$. However, the renormalized field is
  {\em not} divergence--free:
   \[
   \nabla \cdot \left(  \frac{\bmu}{\|\bmu\|} \right)
   = - \frac{a \, b \, b'}{(a^2 + b^2)^{3/2}}.
   \]
  The potentials obtained by solving (\ref{poisson}) are clearly not the same 
  for the original and renormalized fields. The effect of renormalization on 
  the demagnetization field is particularly significant near domain boundaries.
  For example, if the function $b$ appearing in (\ref{renorm_exam}) is a step 
  function, then the divergence of the renormalized magnetic field is a delta 
  function. We thank F.~Reitich, \cite{reitich}, for
  this illustration of the dangers of normalizing vector fields.
   
 \section{Lie group methods for the system $ \dot{\bmu} = \bA(\bmu)\times \bmu$}

  We recall that the general system under consideration is 
  \begin{equation}\label{generalsystem} \frac{\partial}{\partial t}\bmu({\mathbf{x}},t) = {\bf
  A(\bmu(\bx, t))}\,\times\,\bmu({\mathbf{x}},t),\qquad \bmu({\mathbf{x}},0) =  \bmu_0({\mathbf{x}}), \qquad {\mathbf{x}}
  \in {\mathcal B}.\end{equation}
 The initial condition $\bmu_0$ satisfies
  $\|\bmu_0({\mathbf{x}})\|=1$ for all ${\mathbf{x}}\in {\mathcal B} $, i.e. $\bmu_0: {\mathcal B}  \to S^2$, the unit sphere
  in $\Rset^3$. Since $\frac{\partial}{\partial t} \| \bmu({\mathbf{x}},t) \| = 0$ for 
all ${\mathbf{x}}$ and $t$, $\bmu(\ , t): {\mathcal B}  \to S^2$ for all $t$, the rotation group 
$SO(3)$ acts transitively on $S^2$, i.e. any
  point on the sphere can be rotated onto any other point on the sphere; 
  hence there are time dependent curves $\tilde{Q}:{\mathcal B}  \rightarrow SO(3)$
  satisfying
  \begin{equation} \label{rotategen}\bmu({\mathbf{x}}, t) =\tilde{Q}({\mathbf{x}}, t) \bmu_0({\mathbf{x}}),
  \qquad \forall \ \mathbf{x}\in {\mathcal B}. \end{equation}
 We emphasize that these curves are not uniquely determined by (\ref{rotategen}). 

  In the spatially discretized version of the system, we choose 
  $N$ grid points ${\mathbf{x}}_n, n=1, 2,...,N$ in ${\mathcal B}$ and consider 
  only the values of the magnetic field at those grid points. 
  Given a curve $\bmu(\cdot, t): {\mathcal B}  \to S^2$, we define 
 \[ \bM(t):=\left(\bmu({\mathbf{x}}_i,t)\right)_{ i=1}^N \in 
 {\mathcal M} := (S^2)^N.\] 
  Where there is no confusion, we suppress the argument $t$.
  The fully discretized version of the system (\ref{generalsystem}) is
  \begin{equation}\label{discretesystem} \dot \bM = {\bf A(\bM)}\times
  \bM.\end{equation} 
Here and throughout the paper, vector operations such as cross or inner products
 on $(S^2)^N$ or $(\Rset^3)^N$ should be understood as the usual operations in 
$\Rset^3$ performed on each of the $N$ component vectors. The Lie group ${\mathcal G } = (SO(3))^N$ acts
  transitively, but not freely, on ${\mathcal M}$. That is, any point in $\mathcal M$
  can be mapped onto any other point in ${\mathcal M}$ by the action (by component-wise 
  rotations) of ${\mathcal G }$ on ${\mathcal M}$, but the group element accomplishing any such   
  transformation is not unique. The isotropy subgroup of a point $\bM$ in
  ${\mathcal M}$ ( where $\mathcal M$ is the group of transformations fixing $\bM$) is an 
  $N$--dimensional torus. As in (\ref{rotategen}), there are 
  smooth time-dependent curves $Q$ in the group ${\mathcal G }$ satisfying
  \begin{equation} \label{rotatedisc} 
  \bM(t) = Q(t)\bM(0),
  \qquad \mbox{and hence} \qquad
  \dot{\bM}(t) = \dot Q(t)\bM(0).
  \end{equation} 
  
  We can use the Lie algebra 
  \[
  {{\mathfrak G}} = (so(3))^N = \left\{ \mbox{skew symmetric $3 \times 3$ matrices} \right \}^N
  \approx ({\Rset}^3)^N
  \] 
  of ${\mathcal G }$, which is the tangent space to ${\mathcal G }$ at the identity, 
  to put (\ref{rotatedisc}) into a more familiar and convenient form. The
  identification of $so(3)$ with $\Rset^3$ is implemented using the map
   $\mbox{skew}: \Rset^3 \to so(3)$ given by
  \[ 
  \skw{\bxi} :=\left(\begin{array}{cccc}
  0 &-\xi_3& \xi_2\\
  \xi_3&0&-\xi_1\\
  -\xi_2&\xi_1&0\end{array}\right),
  \]
  i.e. $\skw{\bxi} \bx = \bxi \times \bx$ for all 
  $\bx \in \Rset^3$.
  The matrix commutator bracket on $so(3)$ corresponds to the cross product
  on $\Rset^3$ under this identification.

  Given any differentiable curve $Q(t)$ in ${\mathcal G }$, there exists a curve
  $\bxi(t) \in (\Rset^3)^N$ satisfying
  \[
   \dot Q(t) = \skw{\bxi(t)} Q(t), 
  \]
  where the product of $\skw{\bxi}$ and $Q$ is the usual matrix product. 
  Thus the system (\ref{discretesystem}) is equivalent to 
\begin{equation}\label{geom_eqs1}
\dot \bM(t)  =\dot Q(t) \bM_0=\skw{\bxi(t)}Q(t) \bM_0=\bxi(t)\times \bM(t).
\end{equation}
Comparing (\ref{geom_eqs1}) to (\ref{discretesystem}), we see that
\begin{equation}
\label{geom_eqs}
\bxi(t) = \bomega(\bM(t)),
\qquad \mbox{where} \qquad
\bomega(\bM)= {\bf A}(\bM) + \sigma(\bM) \bM
\end{equation}
for an arbitrary scalar function $\sigma: {\mathcal M} \to \Rset$.
The flexibility in the choice of map $\sigma$ arises from the non--freeness
of the action of $SO(3)$ on $S^2$, and thus the action of ${\mathcal G }$ on ${\mathcal M}$; 
distinct ODEs
\[
\dot Q = \skw{\bomega(Q \, m)} Q
\qquad \mbox{and} \qquad
\dot Q = \skw{\widetilde {\bomega}(Q \, m)} Q,
\]
where $\bomega(m) - \widetilde \bomega(m) \in \mbox{span}[m]$ for all $m \in S^2$,
will typically have distinct solution curves in $SO(3)$, but the images in 
$S^2$ of those solution curves under the map $Q \mapsto Q \cdot m_0$ will 
coincide. 

When numerically simulating (\ref{discretesystem}), we want a
time-stepping method that ensures that $\bM_n \in {\mathcal M}$, i.e.
that the norms of the component vectors are identically equal to one.
We can regard (\ref{geom_eqs1}) and (\ref{geom_eqs}) as defining a family of 
ODEs 
\[
\dot Q = \skw{\bomega(Q \, \bM_0)} Q,
\]
parametrized by $\bM_0 \in {\mathcal M}$ and $\sigma: {\mathcal M} \to \Rset$, on 
the group ${\mathcal G }$ and use the techniques developed for geometric
integration on Lie groups to determine approximate discrete solution
curves of these ODE (see \cite{simo1,simo2,MKRK,RKMK} and the references 
therein). When combined with the action of ${\mathcal G }$ on ${\mathcal M}$, these techniques
yield geometric integration schemes for (\ref{generalsystem}) that exactly
preserve the constraint $\bM_n \in {\mathcal M}$, regardless of the step size or the
order of the integrator.

The key idea is the following. Suppose we are given a (right) trivialized 
form $\dot g = \xi(g) g$ of an ODE on a Lie group $G$ and an {\it algorithmic
exponential} \, $\Exp: \mathfrak{g} \to G$ mapping the Lie algebra $\mathfrak{g}$ of $G$ into $G$. 
Then an integrator of order $k$ corresponds to an update of the
form \[ g_{n+1} = \Exp(F(g_{n - p}, \ldots, g_n, \Delta t))g_n,\] 
for some map $F: G^{p + 1} \times \Rset \to \mathfrak{g}$ determined by the 
algorithm and the generator $\xi$. Here we consider only one step 
methods, with $p = 0$. We emphasize that the algorithmic exponential 
need not be the true exponential of the Lie group, or even 
a good approximation to the
true exponential; all that we require is that it map algebra elements
exactly into the group $G$ and that the algorithmic exponential of the
zero vector in $\mathfrak{g}$ equal the identity element 
of $G$. For example, the Cayley transform, given by 
$$\mbox{cay}(\bxi) = (I + \skw{\bxi/2})(I - \skw{\bxi/2})^{-1},$$
is an algorithmic exponential for the rotation group that has long been 
used in computational mechanics. The Cayley transform has long been used to 
implement exact rotations in elasticity and plasticity simulations; see, for 
example, \cite{SVQ,SFox}). More recently, it has been utilized in 
the geometric integration of a wide variety of mechanical systems, including
the LLG equations; see \cite{simo1,simo2,lewis,nigamlewis2000,krishna}, 
and references therein. Algorithms of arbitrarily high order can
be constructed using the Cayley transform, despite the fact that is only
a second order approximation of the matrix exponential of $SO(3)$. 
(Note that the Cayley transform is, in fact, an algorithmic exponential
for any matrix group determined by a quadratic constraint. See, e.g. 
\cite{weyl}.)
For the rotation group $SO(3)$ on $\Rset^3$, both the true matrix exponential
and the Cayley transform can be efficiently evaluated and have frequently 
been used as algorithmic exponentials. The true exponential  takes the form
\[
\exp(\bxi) = I + \frac{\sin \|\bxi\|}{\|\bxi\|} \, \skw{\bxi} 
	+ \frac{1 - \cos \|\bxi\|}{\|\bxi\|^2} \, \skw{\bxi}^2.
\]
The image of a vector $\bx \in \Rset^3$ under the action of $\mbox{cay}(\bxi)$
takes the simple, readily evaluated form
\[
\mbox{cay}(\bxi)\bx
= \bx +\frac{1}{1+\|\bxi/2\|^2}\left[\bxi\times  \bx 
+\half \, \bxi\times(\bxi\times \bx)\right].
\] 
Hence we shall use the Cayley transform (actually, $N$ copies of the 
Cayley transform) as our algorithmic exponential $\Exp: {\mathfrak G} \to {\mathcal G }$.
 We observe that the Cayley transform has
the advantage that the entries of $\mbox{cay}(\bxi)$ are rational functions
of the components of $\bxi$; in particular, no trigonometric functions 
need be evaluated.
  
  To summarize this section, we have rewritten the discrete system \[
  \dot\bM(t) = {\bf A}(\bM)\times \bM(t), \qquad \bM(0)=\bM_0\] as an
  ODE 
\[ \dot Q(t) = \skw{\bomega(Q(t) \bM_0)}Q(t).\] 
on the Lie group ${\mathcal G }= (SO(3))^N$. 
We now need to describe choices of infinitesimal update maps $F$ and
  generators ${\bf A}$ that determine one step numerical updates
  $Q_{n+1}=  \Exp(\bF(Q_n \bM_0, \Delta t))Q_n$ and associated updates
  \[
\bM_{n + 1} = Q_{n + 1} \bM_0 = \Exp(\bF(\bM_n, \Delta t)) \bM_n
  \] with specified properties,
  e.g. a specified order of overall accuracy. The construction of suitable
  updates is the subject of the next section. 
  
  \section{Generator selection --- Expansions, curvature, and partial connections}
\label{gen_crap}

A natural and obvious goal in selection of a numerical scheme is the 
achievement of the highest possible accuracy working within the given 
constraints. However, the prioritization of the constraints (efficiency,
stability, developer effort, preservation of key features of the modeled
system, etc.) can lead to significantly different approaches to the achievement
of this goal and correspondingly different schemes. For the purposes of
this discussion, we shall assume that we are given a family of one step
methods of the form 
\begin{equation} \label{our_kind}
(\bM_{n+1})_j = \mbox{cay}(\bF(\bM_n, \Delta t)_j)(\bM_n)_j
\qquad j = 1, \ldots, N
\end{equation}
for some map $\bF: \mathcal M \times \Rset \to (\Rset^3)^N$. 
(Recall that $\mathcal M :=(S^2)^N$.)

We can optimize accuracy within this family of methods by selecting an 
appropriate isotropy algebra correction; specifically, we shall 
first compute the update generator given some `default' choice of 
generator, and then use that generator to determine an element of the
isotropy algebra of the current state that minimizes the discretization
error for that update. We will derive conditions specifying this choice
of isotropy element for algorithms on $\mathcal M$ utilizing the action of
$\mathcal G$ using a traditional series expansion approach to the 
computation of the discretization error. (See \cite{lo} for
the application of this approach and those discussed below to general
homogeneous manifolds.) Subsequently, we will discuss some more geometric 
and, in some cases, less computationally intensive, approaches to this task. 

The geometric approach most closely 
related to the naive series expansion treatment is the use of geodesic
curvature to characterize the essential information about curves on 
manifolds, e.g.; solution curves of differential equations. A less
directly related approach, but one that coincides with the direct
error minimization approach for the forward Euler method on $\mathcal M$
is the use of a partial moving frame and the associated partial connection 
form to determine a choice of generator. We shall define these constructions
in section \ref{partial_conn}.

\subsection{Algorithms on $S^2$}

To demonstrate the influence of the isotropy algebra on discrete
trajectories, we first study the action of the rotation group
$SO(3)$ on a single sphere $S^2$. The techniques we use in analysing
the case of a single sphere immediately generalize to $\mathcal M$
An autonomous vector field $X$ on $S^2$ satisfies $\left \langle X(m), m\right \rangle = 0$
for all $m \in S^2$; hence there exists a (nonunique) map $A: S^2 \to
\Rset^3$, called a {\em generator} of $X$, satisfying
\[
X(m) = A(m) \times m
\]
on $S^2$. Recall from our earlier discussion that distinct choices
of generator $A$ typically yield distinct discrete trajectories 
when used in a numerical algorithm of the form (\ref{our_kind}).
We will compare the performance
of various schemes with two different choices of generator.
The first choice is the `default' or `natural' one, which is not
assumed to have any particular geometric properties.
For the rigid body system, we will take the body angular velocity
as our default generator. 
For arbitrary vector fields $X$, there need not be a natural choice
of generator; 
our intent here is to make a plausible choice of generator that one 
might make if the issue of isotropy were not taken into account. The 
second choice is the orthogonal generator, i.e. the unique map 
$A_o \colon S^2 \to \Rset^3$ satisfying
\[
X(m) = A_o(m) \times m
\qquad \mbox{and} \,\,
\left \langle A_o(m), m\right \rangle = 0
\]
for all $m \in S^2$. 

A general and direct, but potentially computationally intensive, approach to the choice
of a generator that reduces the discretization error is to compute the lowest 
order nonzero term in the series expansion for the discretization error for
the family of algorithms under consideration, leaving the isotropy algebra
component of the generator as an undetermined parameter, and thus determine
conditions on the isotropy component that minimize the error.
See \cite{lo} for a general treatment of this approach and
section \ref{opt} for the derivation of the optimal generator choice for our
forward Euler algorithm for the LLG system.
 Here we briefly explore some 
alternatives to this approach that have natural geometric interpretations.

We consider an Euler update of the form
\begin{equation}
\label{simple_up}
\widetilde {\mathcal{F}}_{\Delta t}(m) = \Exp(\Delta t\, (A_o(m) + \sigma(m) \, m)),
\end{equation}
where $\sigma$ determines the isotropy (normal component) contribution.
In the special case that the algorithmic exponential is given by a rescaling 
of the usual matrix exponential, e.g. by the Cayley transform, (\ref{simple_up})
satisfies 
\[
\widetilde {\mathcal{F}}_{\Delta t}(m) 
= \exp(\tau(m, \Delta t) \, (A_o(m) + \sigma(m) \, m))
\]
for some rescaling $\tau$ of time. Hence, in 
this case, $\widetilde {\mathcal{F}}_t(m)$ is given by a rigid rotation of $m$ 
about an axis depending only on $m$. Hence the curve
\[
\Gamma_\epsilon(m) = \left\{ \widetilde {\mathcal{F}}_t(m) : |t| \leq \epsilon \right \}
\]
is a segment of a circle in $S^2$. Our goal is to choose $\sigma$
so as to obtain the best circular approximation at $m$ to the true orbit 
segment 
\[
{\mathcal O}_\epsilon(m) = \left\{  {\mathcal{F}}_t(m) : |t| \leq \epsilon \right \}.
\]
If $X(m) \neq 0$, then the optimal circular approximation to 
${\mathcal O}_\epsilon(m)$
at $m$ can be characterized using the geodesic curvature 
\[
k_g(m) = \frac {\left \langle(X \cdot \nabla)X(m), m \times X(m)\right \rangle} {||X(m)||^3}
\]
of ${\mathcal O}_\epsilon(m)$ at $m$. 
The best circular approximation to ${\mathcal O}_\epsilon(m)$ at $m$ 
is tangent to $X(m)$ at $m$ and has geodesic curvature equal to that of
${\mathcal O}_\epsilon(m)$ at $m$. The first condition is clearly satisfied 
for any consistent update. The
geodesic curvature $\tilde k_g(m)$ of $\Gamma_\epsilon(m)$ is easily seen 
to satisfy $|\tilde k_g(m)| = |\cot \phi|$, where $\phi$ is the angle
between $m$ and $A_o(m) + \sigma(m) \, m$ (see, e.g. \cite{doCarmo}, p. 249); thus
$|\tilde k_g(m)| = |\sigma(m)|/||A_o(m)||$. Hence optimal 
orbit capture within the class of updates (\ref{simple_up}) is obtained using
\begin{equation}\label{optimal}
\sigma_{\rm cor}(m) := k_g(m) ||X(m)||.
\end{equation}

If ${\mathcal O}_\epsilon(m)$ is itself a segment of a circle, then $\sigma_{\rm cor}$
yields $\Gamma_\epsilon(m) = {\mathcal O}_\epsilon(m)$. Hence any
torsion--free orbits, e.g. the separatricies of the reduced rigid body
equations, are captured exactly by this version of the Euler method.
Note that the choice $A_o$ is suboptimal for the Euler method unless 
$k_g \equiv 0$ along the tractory of interest, i.e. unless the desired 
trajectory is a great circle.

For higher order methods, the axis of rotation used in the update map $\widetilde {\mathcal{F}}_t$
is typically time dependent and hence the corresponding algorithmic trajectory
segment typically is not circular (i.e. it has nonzero torsion). Hence the
simple argument used in the preceding paragraph cannot be applied. However,
the strategy of curvature--matching can still be followed.
Since a smooth curve on a two dimensional manifold in $\Rset^3$ is determined
up to a time reparametrization by its geodesic curvature, we can determine
the conditions on the generator imposed by the restriction that the geodesic 
curvature of $\widetilde {\mathcal{F}}_t(m)$ match that of ${\mathcal{F}}_t(m)$ to some order. The 
higher order derivatives of the curvature can either be determined analytically
for a given vector field $X$ or numerically approximated using standard
difference schemes. 

More generally, the choice a generator of a vector field on a homogeneous 
manifold can be viewed as a special case of the choice of a partial 
connection form, which generalizes to nonfree actions the classical 
connection form on a principal bundle. A partial connection form is a
Lie algebra--valued one--form with appropriate equivariance properties.
In section \ref{partial_conn} we state the relevant definitions and 
present a family of partial connection forms on open subsets of $S^2$ that 
yield a discretization error--minimizing algorithm
and determine an algorithm that captures orbits to second order
for any dynamical system on a single copy of $S^2$. The interested reader is
referred to \cite{lno} for a more detailed treatment of partial
connection forms and related constructions.

\subsection{Partial connection forms}
\label{partial_conn}

We now briefly discuss a general geometric approach to the selection of 
generators, using a generalization of the connection form on a principal
bundle. For a more detailed treatment of this generalization and proofs 
of the assertions given below, see Lewis et al. [2002].
Let $\mathcal P$ be a principle bundle, that is, a manifold $\mathcal P$ acted on by a Lie group $G$. Let the action of $g\in G$ in $\mathcal P$ be denoted by  $  g \cdot p = \Phi_g(p)= \hat \Phi_p( g),\,\forall  g \in G, p\in\mathcal P.$ Let $\mathfrak{g}:= T_eG$ denote the algebra of $G$ and let $\mathfrak{g}\cdot p:=\left\{ T_e\hat\Phi_p\cdot \xi :\,\xi\in\mathfrak{g}\,\right\}$.
Recall that a {\em connection} on a principal bundle ${\mathcal P}$ is a distribution 
$\Gamma$ satisfying
\[
T_p {\mathcal P}= \mathfrak{g} \cdot p \oplus \Gamma_p 
\qquad \mbox{and} \,\,
T_p \Phi_g \cdot \Gamma_p = \Gamma_{g \cdot p},
\]
for all $p \in {\mathcal P}$ and $g \in G$. Specification of a 
connection $\Gamma$ is equivalent to specification of an equivariant 
$\mathfrak{g}$--valued one--form $\alpha$, called the {\em connection form}, satisfying
\[
\alpha \circ T_e \hat {\Phi}_p = \mbox{id}, \qquad \qquad 
i.e. \quad \alpha(p)(\xi_{\mathcal P}(p)) = \xi \qquad \mbox{for all $\xi \in \mathfrak{g}$},
\]
for all $p \in {\mathcal P}$. By equivariance we mean that $\alpha \circ T \Phi_g
= \mbox{Ad}_g \circ \alpha$ for all $g \in G$. The connection $\Gamma$ and
connection form $\alpha$ are related by the condition
$\mbox{ker}[\alpha(p)] = \Gamma_p$ for all $p \in {\mathcal P}$. 
(See, e.g., \cite{kn} for a detailed presentation of the 
properties of connections and connection forms.)

The equivariance properties of connections and connection forms typically
cannot be preserved in the context of nonfree actions, hence we relax these 
conditions, requiring only equivariance with respect to specified 
representatives of the isotropy equivalence classes.
A map $\beta: G \times {\mathcal M}\to G$ is a {\em slip map} if $\beta(g, m) \cdot m 
= g \cdot m$ for all $g \in G$ and $m \in {\mathcal M}$. 
A (singular) distribution $\aleph$ assigning a complement 
$\aleph_m$ to $\mathfrak{g} \cdot m$ in $T_m {\mathcal M}$ to each point $m \in {\mathcal M}$ is a 
{\em partial connection} if there is a slip map $\beta$ satisfying
\[
T_m \Phi_{\beta(g, m)} \cdot \aleph_m = \aleph_{g \cdot m}
\]
for all $g \in G$ and $m \in {\mathcal M}$. 
A {\em partial connection form} with slip map $\beta$ is a 
$\mathfrak{g}$--valued one--form $\alpha$ on ${\mathcal M}$ satisfying
\[
\alpha(m)\left({\eta}_{\mathcal M}(m)\right) = \eta \quad \mbox{mod} \quad \mathfrak{g}_m, \qquad
\mbox{i.e.} \qquad T_e \hat \Phi_m (\alpha \circ T_e \hat \Phi_m - \mbox{id}) = 0
\]
and
\[
\Phi_{\beta(g, m)}^*\alpha(m) = \mbox{Ad}_{\beta(g, m)}\alpha(m) 
	\quad \mbox{mod} \quad \mathfrak{g}_{g \cdot m}
\]
for all $g \in G$ and $m \in {\mathcal M}$.

One natural source of partial connections is a generalization of a
{\em moving frame}, in the modern sense introduced by Fels and Olver (\cite{FOmcI}, 
\cite{FOmcII}), i.e. a smooth equivariant map $\rho: {\mathcal P}\to G$ on a manifold  
${\mathcal P}$ with a free group action. 
Recall that a principal bundle ${\mathcal P}$ is trivial if there exists a global
section, i.e. a smooth map $\Sigma: {\mathcal M}\to {\mathcal P}$ from the base manifold ${\mathcal M}$ 
into${\mathcal P}$ such that each group orbit $G \cdot p$ in ${\mathcal P}$ intersects 
$\Sigma({\mathcal M})$ exactly once and the projection $\pi: {\mathcal P}\to {\mathcal M}$ satisfies
$\pi \circ \Sigma = \mbox{id}$. This condition corresponds to the existence of 
a moving frame. Specifically, a global section $\Sigma$ and associated moving 
frame $\rho$ are related by the equality
\[
\Sigma(\pi(p)) = \rho(p)^{-1} \cdot p 
\]
for all $p \in {\mathcal P}$. A global section determines a flat connection, namely
the connection that assigns to any point $\Sigma(m)$ the subspace 
$T_m \Sigma \cdot T_m {\mathcal M}$. 
If we introduce the notation $\tilde D \gamma: T {\mathcal P}\to \mathfrak{g}$ to denote the
right trivialization of the linearization of a map $\gamma: {\mathcal P}\to G$, i.e.
\[
\tilde D \gamma(\delta p) := T_p (R_{\gamma(p)^{-1}} \circ \gamma) \, \delta p
\]
for any $p \in {\mathcal P}$ and $\delta p \in T_p {\mathcal P}$, then $\tilde D \rho$ is the
connection form of the connection determined by $\Sigma$. 

The rotation group $SO(3)$ acts transitively on $S^2$ and freely and 
transitively on the unit tangent bundle $U(S^2) = \left\{ u \in T S^2 : 
\|u\| = 1 \right \}$. The map $\rho: U(S^2) \to SO(3)$ taking $u \in U_m S^2$
to the orthogonal matrix with columns $(m, u, m \times u)$ is a (left) 
moving frame with associated connection form
\[
\tilde D \rho(\delta u) = m \times \delta m + \left \langle u \times \delta u, m\right \rangle m,
\]
where $\delta u\in T_u U(S^2)$, with $m = \pi(u)$ and $\delta m = T_u \pi \, \delta u$.
(Here $\pi: U(S^2) \to S^2$ denotes the canonical projection.)
Note that we will regard $u$ both as a tangent vector
to the sphere at $m$ and as a unit vector in $\Rset^3$.

Moving frames can be extended to manifolds with nonfree actions as follows:
A (smooth) map $\phi: {\mathcal M}\to G$ is a (left) {\em partial moving frame} if 
\begin{equation}\label{partial moving}
\phi_g(m) := \phi(g \cdot m) (\phi(m))^{-1} 
\end{equation}
satisfies
\[
\phi_g(m) \cdot m = g \cdot m
\]
for all $g \in G$ and $m \in {\mathcal M}$. A {\em partial moving frame} on a 
submanifold ${\mathcal S}$ of a manifold ${\mathcal M}$ with a $G$ action is a map 
$\phi: {\mathcal S} \to G$ satisfying (\ref{partial moving}) for any $m \in {\mathcal S}$ and 
any $g \in G$ such that $g \cdot m \in {\mathcal S}$.
The trivialized linearization $\tilde D \phi$ of a partial moving frame 
$\phi: {\mathcal M}\to G$ is a partial connection form, with associated slip map
$\beta(g, m) = \phi_g(m)$.
We refer to the trivialized linearization $\tilde D \phi$ as the 
partial connection form associated to the partial moving frame $\phi$. 

If a group $G$ acts transitively on a manifold ${\mathcal M}$, then every group 
orbit is equal to the entire manifold. In this situation, a (partial)
connection form $\alpha$ assigns to each tangent vector a generator of 
that vector, so that
\[
(\alpha(m)(\delta m))_{{\mathcal M}}(m) = \delta m 
\]
for all $\delta m \in T_m {\mathcal M}$ and all $m \in {\mathcal M}$. In particular, given a 
vector field $X$ on ${\mathcal M}$, the map $\omega := \iota_X \alpha: {\mathcal M}\to \mathfrak{g}$, 
i.e. $\omega(m) = \alpha(m)(X(m))$, satisfies
\[
\omega(m)_{{\mathcal M}}(m) = X(m)
\]
for all $m \in {\mathcal M}$. Hence (partial) connection forms can be used to construct
geometric integration schemes on manifolds with transitive actions. We are
currently investigating the role of geometrically motivated choices of
partial connection forms in the design of efficient geometric integration
algorithms. As can be seen in the example discussed below, simple, natural
choices of partial connection forms can lead to significant improvement
in numerical performance.

To illustrate the somewhat abstract geometric constructions described above, we 
now present a moving frame associated to the action of the rotation group
$SO(3)$ on the unit tangent bundle $U(S^2)$ of the sphere $S^2$ and an
associated family of partial moving frames on $S^2$. The partial connection
form (which, in this case, is simply a map from the sphere to $\Rset^3$)
of one of these partial moving frames yields the discretization 
error--minimizing generators used in the versions of the forward Euler method 
described in sections \ref{opt} and \ref{rig_body}.
As we shall see, this partial connection form and the associated 
generators can be derived without the use of the expansion of the 
discretization error.

Any unit vector field $Y$ on a submanifold ${\mathcal M}$ of $S^2$ 
determines a partial moving frame \ $\phi = \rho \circ Y$ on ${\mathcal M}$, with 
partial connection form 
\begin{equation}\label{fuzzconn_sphere}
\tilde D \phi(\delta m) = m \times \delta m + \left \langle Y(m) \times (D Y(m) \cdot \delta m), m\right \rangle m.
\end{equation}
The map $\phi_g$ associated to $g \in SO(3)$ is
$\phi_g(m) = g \, \exp(\theta(g, m) \, m)$, where
$\theta(g, m)$ denotes the angle between $g^{-1} Y(g \, m)$ and $Y(m)$.
Given an ODE $\dot m = X(m)$, we can set
$Y(m) = X(m)/\left\Vert{X(m)}\right\Vert$ on some set ${\mathcal M}\subset S^2$ containing no
equilibria (i.e. zeroes of $X$); in this case, (\ref{fuzzconn_sphere}) takes the form
\begin{equation}\label{fuzzconn_sphere2} \displaystyle
\tilde D \phi(m)(\delta m) = m \times \delta m + 
\frac{\left \langle (\delta m \cdot \nabla) X, m \times X(m)\right \rangle}{\left\Vert{X(m)}\right\Vert^2} \, m.
\end{equation}
In particular,
\begin{equation}\label{fuzzconn_sphere3}
\tilde D \phi(m)(X(m)) = m \times X(m) + k_g(m) \,\left\Vert{X(m)}\right\Vert m,
\end{equation}
where $k_g(m)$ denotes the geodesic curvature of the curve $m(t)$
in $S^2$. 

The partial connection form (\ref{fuzzconn_sphere}) can be used to select the
isotropy correction map $\sigma$ used in (\ref{simple_up}).
Following (\ref{fuzzconn_sphere3}), we set 
\begin{equation}\label{sig_choice}
\sigma(m)(x) := \left\{ \begin{array}{ll}
k_g(m)(x) \, \Vert{X(m)(x)}\Vert \qquad & X(m)(x) \neq 0 \\
0 & X(m)(x) = 0
\end{array} \, . \right .
\end{equation}
In our numerical implementation (\ref{sig_choice}), we approximate $k_g(m)$
using the identity
\[
k_g(m) = \frac{\left \langle\ddot m, m \times \dot m\right \rangle}{||\dot m||^3}
= -  \frac{\left \langle\dot \omega, \omega \times m\right \rangle}{||\omega||^3}
\]
for a curve $m(t)$ in $S^2$ with nonzero velocity $\dot m = \omega \times m$,
where $\omega$ is orthogonal to $m$, and replacing $\dot m$ and
$\ddot m$ with finite difference approximations.


   \section{ A collection of geometric integrators} \label{sectionalgo}
 Recall that an update of the form
$$Q_{n + 1} = \Exp(\bF(Q_n \bM_0, \Delta t)) Q_n,$$
where $\Exp: {\mathfrak G} \to {\mathcal G}$ satisfies $(\Exp(\bxi))_j = 
\mbox{cay}(\bxi_j)$, determines a one 
step method on the Lie group ${\mathcal G }$, where $\bF: {\mathcal G } \times \Rset \to {\mathfrak G}$ 
is determined by the generator and the selected scheme, 
and $\Delta t$ denotes the time step. Given an ODE 
$ \dot \bM(t) = {\bf A}(\bM)\times \bM $
on ${\mathcal M}$, we will construct updates of the form
\begin{equation}
\label{gen_update}
\bM_{n+1}=\Exp(\bF(\bM_n, \Delta t))\cdot \bM_n,
\end{equation}
It is our objective to identify several classes of infinitesimal update maps $F$, 
leading to algorithms of first, second, and fourth order.

\subsection{Discretization error reduction through choice of $\sigma$}
\label{opt}
  
As was previously discussed, the normal component of the generator ${\bf A}$ does not influence 
the solution curves of the original ODE. Thus, if we have a numerical 
algorithm of order $n$, this component does not affect the solution up
to order $n$. However, it typically does appear in the higher order terms 
of the approximation, and theoretically a suitable choice of this component 
will reduce the discretization error.
 For the forward Euler scheme an optimal choice of $\sigma$, in
 the sense that this choice minimizes the discretization error, also has a 
 natural geometric interpretation. Here we derive this map $\sigma$ using
 a direct discretization error calculation; in the following section we shall
 discuss various geometric considerations that can be used in the selection
 of the generator to be used in a Lie group integration scheme. 
  
We consider consistent algorithms using standard methods on the tangential 
component of $\mathbf{\xi}$. The normal component is treated as a function of the 
tangential one; we shall see that 
a component of the local discretization error at second order 
can be eliminated by a suitable choice of the normal component 
of the lowest order term in $\mathbf{\xi}$. For the sake of simplicity, we consider
here only the lowest order case, in which the discretization error of a first 
order method is reduced by an appropriate selection of $\sigma$. This is a
 particular example of a more general result covering a large class of 
 manifolds and higher order algorithms. (See \cite{lo,lno}.)
 Work is in progress (Lewis, Nigam, Olver) to possibly extend these or
 related results to an even larger class of systems, including the full 
 discretized LLG system.

We begin by examining the flow ${\mathcal F}_t$ of the ODE
$\dot\bM=\mathbf{A}(\bM) \times \bM$ on ${\mathcal M}$. 
This flow satisfies 
\[
{\mathcal F}_{\Delta t}(\bM)=\bM + \Delta t \, \mathbf{A}\times \bM +\frac{\Delta t^2}{2}  
\left(\mathbf{A}\times (\mathbf{A}\times \bM)+\dot\mathbf{A}\times \bM\right)
+ {\mathcal O}(\Delta t^3).
\]

If the algorithmic update 
$\widetilde {\mathcal F}_{\Delta t}:{\mathcal M}\rightarrow {\mathcal M}$ is given by 
\[ 
\widetilde {\mathcal F}_{\Delta t}(\bM):= 
\Exp(\bF(\bM,\Delta t)) \cdot \bM
\] 
for some map $\bF(\bM,\Delta t):=\mathop{\sum}\limits_{j=1}^{\infty}
\frac{\Delta t^j}{j!}\bxi_j(\bM)$ and $\Exp: \Rset^3 \approx so(3) \to SO(3)$ agrees 
with the exponential map to second order (e.g. $\Exp$ is the Cayley transform), then
\begin{eqnarray*}
\widetilde {\mathcal F}_{\Delta t}(\bM)
&=& \left(I+\Delta t \, \skw{{\bf A}} +\half\Delta t^2 \skw{{\bf A}}^2 
+ {\mathcal O}(\Delta t^3)\right)\bM\\
&=&\left(I+\Delta t \, \skw{\bxi_1} 
+ \half\Delta t^2 \left (\skw{\bxi_2} + \skw{\bxi_1}^2 \right ) 
	+ {\mathcal O}(\Delta t^3) \right)\bM\\\
&=& \bM+\Delta t \, \bxi_1 \times \bM +\frac{\Delta t^2}{2}\left(\bxi_1\times(\bxi_1\times \bM)+\mathbf{\bxi}_2\times \bM\right)
	+ {\mathcal O}(\Delta t^3). 
\end{eqnarray*}

We now derive conditions on the terms $\bxi_1$, $\bxi_2$,\ldots in the series
expansion of $F$ yielding algorithms of increasingly high order.
The consistency condition for $\widetilde {\mathcal F}_{\Delta t}$ is 
${\mathbb P}_{\bM} (\bxi_1 - \mathbf{A}) = 0$, 
where ${\mathbb P}_{\bM}$ denotes component-wise projection onto the orthogonal 
complements of the component vectors of $\bM$, i.e. 
$\left \langle({\mathbb P}_{\bM} \xi)_j, \bM_j\right \rangle = 0$, $j = 1, \ldots, N$.  
If $\widetilde {\mathcal F}_{\Delta t}$ is consistent, then, setting 
$\sigma_1 := \langle \bxi_1, \bM\rangle$, the local discretization error is 
\begin{eqnarray*} 
\frac{\widetilde {\mathcal{F}}_{\Delta t}(\bM) - {\mathcal F}_{\Delta t}(\bM)}{\Delta t}
&=&\frac{\Delta t}{2}\left( -\sigma_1\bM\times 
(\bxi_1\times\bM)+(\bxi_2-\dot{\mathbf A})\times \bM)\right) + {\mathcal O}(\Delta t^2) \\
&=&\frac{\Delta t}{2}\left(-\sigma_1{\mathbf A} 
	+(\bxi_2-\dot{\mathbf A})\times \bM)\right) + {\mathcal O}(\Delta t^2).
\end{eqnarray*}
The algorithm is thus second--order accurate iff
\begin{equation}
\label{second_order}
\langle\bxi_2-\dot {\mathbf A},{\mathbf A}\rangle=0 \qquad \mbox{and} \qquad
\sigma_1 \langle {\mathbf A},{\mathbf A}\rangle
=\langle (\bxi_2-\dot{\mathbf A})\times\bM,{\mathbf A}\rangle.
\end{equation}
In our geometric version of the forward Euler method with $\bF(\bM,\Delta t)=\mathbf{A}(\bM)$, $\bxi_j = 0$ for
$j > 1$; thus this method will not be second order. However, we are 
free to choose $\sigma_1$ so as to satisfy the second equality
in (\ref{second_order}), e.g.
\begin{equation}\label{movingsigma} \displaystyle
\sigma(\bM, \Delta t) 
	= \frac{\left \langle\Delta {\mathbf A}(\bM, \Delta t), {\mathbf A}(\bM) \times \bM\right \rangle}{\left\Vert{{\mathbf A}(\bM)}\right\Vert^2} + {\mathcal O}(\Delta t^2),
\end{equation}
where $\Delta {\mathbf A}(\bM, \Delta t)$ is some first order approximation to
$\dot {\mathbf A}(\bM)$ (e.g., a discrete difference approximation), yielding a 
discretization error--minimizing member of the family of algorithms with $\bF(\bM,\Delta t)={\mathbf A}(\bM)+\sigma(\bM,\Delta t)\bM.$

Analogous expansions can be used to minimize the discretization error of higher
order methods. However, the symbolic calculation of such expansions for high order 
schemes is, at present, relatively laborious and does not seem tractable for systems 
such as the LLG equations, in which the generator $\mathbf{A}$ is determined in part (the
demagnetization field) by a nonlinear PDE. 

\subsection{First order methods} 

Using (\ref{gen_update}), we now define geometric one-step methods that are natural 
analogs of the standard explicit and implicit Euler methods: 
\begin{equation}
\label{method1}
\bF(\bM, \Delta t) = \left \{ \begin{array}{ll} {\bf A}(\bM)+\sigma(\bM,\Delta t) \bM \qquad \quad 
	& \mbox{forward Euler},\smallskip \\
\tilde {\bf A}(\bM, \Delta t)
  + \tilde \sigma(\bM,\Delta t) \bM & \mbox{implicit Euler},
\end{array} \right . 
\end{equation}
where $\tilde{\bf A}(\bM,\Delta t)$ denotes the solution of the implicit equation 
$\xi = {\bf A}(\Exp(\Delta t \, \xi)\bM)$ and the scalar functions $\sigma$
and $\tilde \sigma$ are as yet unspecified. 

The numerical results presented in sections \ref{gen_crap} and \ref{num_results}
illustrate the effect of the parameter $\sigma$ 
on the discrete trajectories determined by the forward Euler algorithm when applied to
rigid body dynamics and the LLG micromagnetism model. We shall see that in the rigid
body system, $\sigma$ satisfying (\ref{movingsigma}) yields second order accuracy in 
energy tracking, and thus second order orbit capture for this conservative system. 
In the micromagnetics simulations, where
damping plays a crucial role in the long term dynamics, large values of
$\sigma$ cause the trajectories to sharply diverge from those of the
ordinary forward Euler; however, the final state is the same. A
closer look at the LLG equation shows that a larger value of $\sigma$
corresponds to the inclusion of more precession in the trajectory.
These numerical results clearly show that different choices of the 
parameter $\sigma$ lead to significantly different numerical trajectories 
and thus motivate the search for an ``optimal''
value of $\sigma$. In section \ref{gen_crap} we describe a general 
geometric approach to selecting values for $\sigma$; in section
\ref{opt} we show that when used with the forward Euler method,
this choice of $\sigma$ minimizes the discretization error.
(See \cite{lo} for a description of this
approach for more general manifolds.)
 
 \subsection{Second order methods}\label{second_ordersection}
  
  We consider four second order methods modeled on the classic Heun (RK2)
algorithm. 
  
In the first, we use the `default' generator $\mathbf{A}$ in the Heun method, i.e.
\[
\bF_{\rm def}^{RK2}(\bM, \Delta t) := \half \left(  \mathbf{A}(\cay{\Delta t\, \mathbf{A}(\bM)} \bM)
	+ \mathbf{A}(\bM) \right).
\]
The second method, $\bF_{\rm orth}^{RK2}$, is entirely analogous, but with 
$\mathbf{A}$ replaced with the orthogonal generator $\mathbf{A}_o$, where 
  \[ \mathbf{A}_o(\bM) :={\bf A}(\bM)- (\bM\cdot {\bf A}(\bM))\bM,\] 
  and hence $\left \langle(\mathbf{A}_o(\bM))_j, \bM_j\right \rangle = 0$, $j = 1, \ldots, N$.
  
  In the third method, the infinitesimal
rotation determined by applying the Heun method to the default generator $\mathbf{A}$
is modified by addition of an appropriate isotropy element to yield a higher
order of orbit capture; specifically,
\[
\bF^{RK2}_{\rm dcor}(\bM, \Delta t) := \bF_{\rm def}^{RK2}(\bM, \Delta t) 
	+ \Delta t^2 \, \sigma_{\rm def}(\bM) \bM 
\]
The fourth algorithm is analogous, but with $\bF_{\rm def}^{RK2}$ replaced by 
$\bF_{\rm orth}^{RK2}$ and $\sigma_{\rm def}$ replaced by an appropriate function $\sigma_{\rm orth}$. 
(The function $\sigma_{\rm orth}$ is a rational function in
$m$ and the components of $\mathbf{A}$, but is significantly more
complicated than $\sigma_{\rm def}$.) Note that isotropy plays a role both in the choice of
the generator and in the selection of a correction term.

  \subsection{ Fourth order methods}
  In this subsection, we describe two families of fourth order Lie group 
 integrators on $\mathcal M$. We emphasize that these algorithms map a 
point $\bM_n \in \mathcal M$ {\em exactly} into $\mathcal M$; they are fourth
order accurate in the sense that they approximate the true trajectories
within $\mathcal M$ to fourth order.
A direct implementation of the classic fourth--order Runge--Kutta method
to a vector field on $\mathcal M$ fails to maps exactly into the manifold
$\mathcal M$, while application of the classical RK4 method to the generator of the
flow, followed by application of  the exponential map and the group action necessarily
yields an update in $\mathcal M$, but typically does not give a fourth order 
approximation of the true flow. The generator must be modified to account
for the trivialization of the tangent bundle of the group; this 
modification can be implemented either before or after the stages of the
Runge--Kutta method are computed and averaged.
  The first method we use is the RKMK4 method,
a Runge--Kutta style method due to Munthe--Kaas (\cite{MKRK}, \cite{RKMK})
 in which each stage of a traditional RK4 method is modified so that
 the resulting generator, followed by (algorithmic) exponentiation and
application of the group action to the manifold, yields a fourth order method.
The second method utilizes a series expansion of the generator along the true 
flow, followed by a single modification to account for the trivialization of the
tangent bundle of the group ${\mathcal G }$, again followed by (algorithmic) 
exponentiation and application of the group action.
  
We implemented the
  RKMK4 method \cite{MKRK}, using the Cayley transform 
  rather than the true matrix exponential.   
  To implement a Lie group integrator for (\ref{discretesystem}) using the Cayley transform, we 
  make use of the fact that, for sufficiently small
  $t$, there is a function $\blf: \Rset \to \Rset^3$ satisfying
  \begin{equation}\label{use_cay}
  \bM(t)= \cay{\blf(t)} \bM(0).
  \end{equation}
  Differentiating (\ref{use_cay}) with respect to $t$, we obtain 
  \begin{eqnarray*}
  \dot \bM(t)
  &=& \mathrm{dcay}_{\blf}(\blf'(t)) \times \cay{\blf(t)}\bM(0) \\
  &=& \mathrm{dcay}_{\blf}(\blf'(t)) \times \bM(t) \\
  &=& \bA(\bM(t)) \times \bM(t),
  \end{eqnarray*}   where the map $\mathrm{dcay}_{\blf} 
  = \tilde D{\mbox{cay}}(\blf): \Rset^3 \to \Rset^3$
  is the right trivialization of the tangent map of the Cayley transform.
  
  Hence $\blf'$ and $\bA$ are related by
  \[
  \mathrm{dcay}_{\blf}(\blf'(t)) \times \bM(t) = \bA(\bM(t)) \times \bM(t),
  \]
  which is equivalent to
  \begin{equation} \label{big} 
  \blf'(t) = \mathrm{dcay}_{\blf}^{-1}(\bA(\bM(t)) + \sigma(t) \bM(t))
  \end{equation}
  for some function $\sigma$. The initial condition
  for (\ref{use_cay}) is $\blf(0)= 0$.
  
  The map $\mathrm{dcay}_{\blf}: \Rset^3 \to \Rset^3$ satisfies
  \[
  \mathrm{dcay}_{\blf}:= \frac{1}{1+\|\half \blf\|^2}(I + \half \, \skw{\blf})
  \]
and
\[
\mathrm{dcay}_{\blf}^{-1} = I -\half \, \skw{\blf} +\frac{1}{4}\blf \, \blf^T.
\]
 The Cayley version of the RKMK4 method is essentially the conventional RK4
  applied to (\ref{big}): having found $\bM_n$ at time $t_n$, we construct
  the update \[ \bM_{n+1}=\bM(t_n+\Delta t)= \cay{\bF^{RK4}(\bM_n, \Delta t)} \bM_n, \] where
  \[ \bF^{RK4}(\bM_n, \Delta t) = \frac{1}{6}(\bF_1+2\bF_2+2\bF_3+\bF_4)\] and  
\[
\begin{array}{ll}
  \bA_1= \Delta t\, {\mathbf{A}(\bM_n)}, & \bF_1= \mathrm{dcay}_{\bf{0}}^{-1}(\bA_1)
	\bigskip \\
  \bA_2= \Delta t\, {\mathbf{A}(\cay{\half \bA_1}\bM_n)}, \qquad & \bF_2=
  \mathrm{dcay}_{\half \bA_1}^{-1}(\bA_2)
	\medskip \\
   \bA_3= \Delta t\, {\mathbf{A}(\cay{\half \bA_2}\bM_n)}, & \bF_3=
  \mathrm{dcay}_{\half \bA_2}^{-1}(\bA_3)
	\medskip \\
   \bA_4= \Delta t\, {\mathbf{A}(\cay{\bA_3}\bM_n)}, & \bF_4=
  \mathrm{dcay}_{\bA_3}^{-1}(\bA_4)
\end{array} .
 \]
For more details on this method, see \cite{MKRK}.

The second method involves a series expansion with respect to time of the 
generator $\bA({\mathcal F}_t(\bM))$, where ${\mathcal F}_t$ denotes the exact flow at 
time $t$. Given the generator $\bA$, this expansion is computed by iteratively 
symbolically differentiating $\bA({\mathcal F}_t(\bM))$ and then substituting 
$\bA(\bM) \times \bM$ for $\dot \bM$. The third order approximation to 
$\bA({\mathcal F}_t(\bM))$ is then modified to take into account the trivialization of
the tangent bundle of $SO(3)$ and the action of $SO(3)$ on $S^2$. The specific
expressions for this modification for the rigid body equations on the sphere
are given in \S \ref{rig_body}. Given the implicit and highly nonlinear nature 
of the LLG equations, symbolic calculation of the derivatives of $\bA$ for 
this system seemed excessively complicated; hence we did not 
implement this algorithm for the LLG system.

\section{An example: the rigid body flow on a sphere}
\label{rig_body}

We now apply the results outlined above to a simple and familiar system,
the reduced rigid body equations on the sphere. Given a positive definite 
symmetric  three by three matrix $\mathbb{I}$, define the vector field 
\begin{equation}\label{rigid_body}
X(m) = m \times \mathbb{I}^{-1} m
\end{equation}
on $S^2$. This is a Hamiltonian system with respect to the 
Kostant--Kirillov-Souraiu symplectic structure 
\[
\Omega(m)(\xi \times m, \eta \times m) = \langle m, \xi \times \eta \rangle
\]
and Hamiltonian 
\begin{equation}\label{rb_Ham}
H(m) = \half \langle m, \mathbb{I}^{-1} m \rangle.
\end{equation}
The system (\ref{rigid_body}) is the symplectic reduction of the free rigid 
body equations on $T^* SO(3)$; more concretely, it is the restriction of
Euler's equation for the body angular momentum to the unit sphere. (Since 
the norm of the body momentum is preserved by the dynamics of Euler's equation,
all spheres centered at the origin are invariant submanifolds.)
The conservative nature of this system makes it particularly easy to measure
the error in orbit capture; if the body is triaxial, i.e. the eigenvalues 
$I_1$, $I_2$, $I_3$ of 
the inertia tensor $\mathbb{I}$ are distinct, the level sets of the Hamiltonian 
(\ref{rb_Ham}) exactly determine the orbits of the system. Thus in this 
situation the error in the orbit is a function of the fluctuation in the energy.
As the numerical results given in tables 1--4 demonstrate, geometric integration
techniques yield efficient, accurate orbit capture for the reduced free rigid 
body, with good performance even for very large time steps. Note that the 
same randomly generated initial conditions and inertia tensors are used in 
all of the numerical simulations.

If the rigid body is axisymmetric, then all true trajectories consist 
either of equilibria (the `poles' and the `equator') or of steady rotations 
in the plane of symmetry. Note that in this situation, even an exactly 
energy--preserving scheme may allow drift across the family of one--point
orbits along the equator. However, all of the methods considered here 
detect equilibria as such. Thus even in the axisymmetric case, we can use
the energy to monitor orbit capture. We shall see that for some of the 
algorithms considered here, there are significant differences in performance
on triaxial and axisymmetric bodies. Symmetries play a crucial role in 
algorithm design and analysis; see, e.g. \cite{lewis}. However, we shall
not explore those issues in any detail here.

\subsection{ Euler methods for rigid body dynamics}

We now consider implementations of the families of algorithms described in
\S \ref{sectionalgo} for the rigid body equations. We take as our default generator
the body angular velocity $A(m) = \mathbb{I}^{-1} m$.
We first consider three first order methods, with infinitesimal updates
\begin{itemize}
\item
$F^{\rm Eul}_{\rm def}(m) = \mathbb{I}^{-1} m$ \medskip

\item
$F^{\rm Eul}_{\rm orth}(m) 
= \mathbb{I}^{-1} m - \langle m, \mathbb{I}^{-1} m \rangle \, m 
= A(m) - 2 \, H(m) \, m$ \medskip

\item
${\displaystyle F^{\rm Eul}_{\rm cor}(m) = \mathbb{I}^{-1} m - \frac {\langle X(m), 
\mathbb{I}^{-1} X(m) \rangle}{||X(m)||^2} \, m
= \mathbb{I}^{-1} m + \frac {\tau(u(m))}{\tau(\mathbb{I} u(m))} \, m}$,

\end{itemize}
where $\tau: \Rset^3 \to \Rset$ and ${u}: S^2 \to \Rset^3$ are given 
with respect to an eigenbasis of $\mathbb{I}$ by
\[
\tau(\bx) = x_1 + x_2 + x_3 
\qquad \mbox{and} \,\,
u(m)_i := (I_j - I_k)^2 I_i m_j m_k
\]
for any cyclic permutation $(i, j, k)$ of $(1, 2, 3)$. 

\newcommand{\leavespace}{\phantom{{\LARGE I}}}
\begin{table}[h]\label{plot_errors}
\begin{center}
\begin{tabular}{| l | c | c | c |} \hline
\leavespace
& $F^{\rm Eul}_{\rm def}$ & $F^{\rm Eul}_{\rm orth}$ & $F^{\rm Eul}_{\rm cor}$ \\ \hline 
\leavespace
Triaxial  & $6.37 \ 10^{-2}$ & $2.23 \ 10^{-2}$ & $4.60 \ 10^{-6}$  \\
\leavespace
Axisymmetric & $2.46 \ 10^{-1}$ & $1.46 \ 10^{-1}$ & $7.38 \ 10^{-14}$  \\
\hline
\end{tabular}
\end{center}
\caption{Maximum energy error over the trajectories given in figure 
	\ref{traj_plots}.}
\end{table}

\setlength{\unitlength}{1in}
\begin{figure}
\label{traj_plots}
\vspace{60mm}
\begin{picture}(6.75,4.5)
\put(0,2.25){\resizebox{2.25in}{2.25in}
	{\includegraphics*[.5in,.5in][3.5in,3.5in]{spheres1n.pdf}}}
\put(2.25,2.25){\resizebox{2.25in}{2.25in}
	{\includegraphics*[.5in,.5in][3.5in,3.5in]{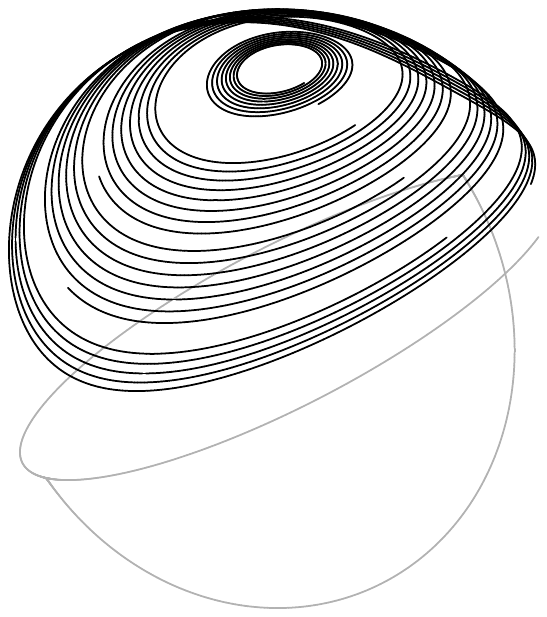}}}
\put(4.5,2.25){\resizebox{2.25in}{2.25in}
	{\includegraphics*[.5in,.5in][3.5in,3.5in]{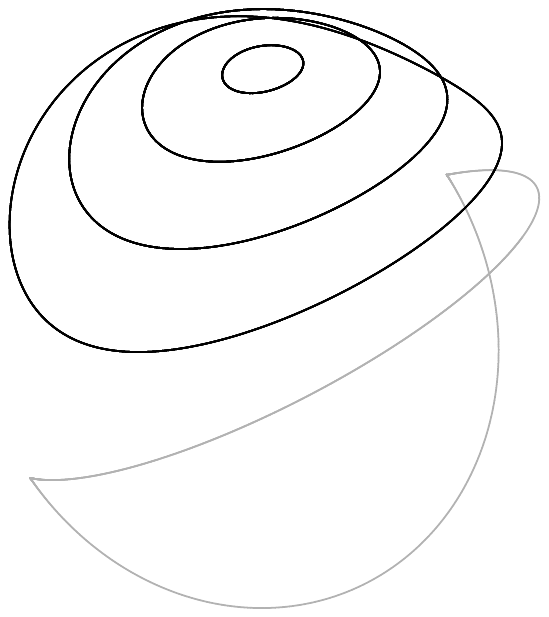}}}
\put(0,0){\resizebox{2.25in}{2.25in}
	{\includegraphics*[0.25in,0.25in][3.75in,3.75in]{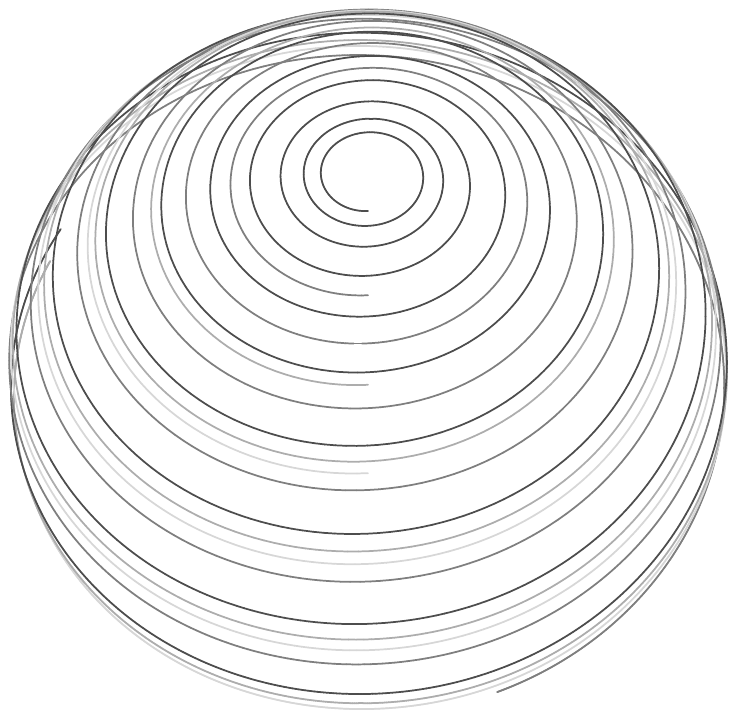}}}
\put(2.25,0){\resizebox{2.25in}{2.25in}
	{\includegraphics*[0.25in,0.25in][3.75in,3.75in]{symspherep1.pdf}}}
\put(4.5,0){\resizebox{2.25in}{2.25in}
	{\includegraphics*[0.25in,0.25in][3.75in,3.75in]{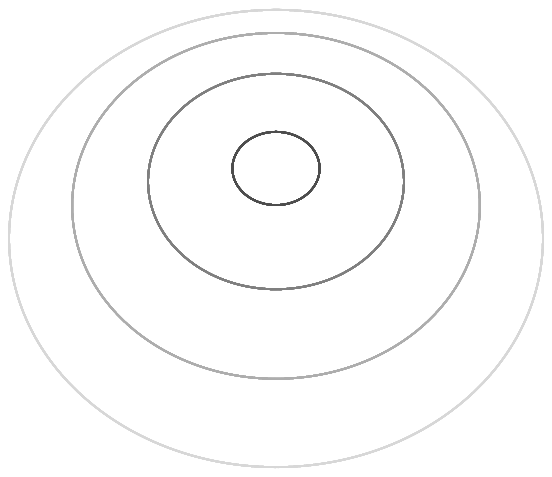}}}
\end{picture}
\caption{Sample trajectories computed over the interval [0, 200] using the
	time step $\Delta t = 0.1$ and, left to right, the first order 
	infinitesimal updates $F^{\rm Eul}_{\rm def}$, 
	$F^{\rm Eul}_{\rm orth}$, and $F^{\rm Eul}_{\rm cor}$.
	The upper row is computed using the inertia tensor of a trixial 
	rigid body, while the lower row is computed for an axisymmetric rigid body.}
\end{figure}

\begin{table}[h]
\begin{center}
\begin{tabular}{| r | c | c | c | c | c | c |} \hline
\leavespace
& \multicolumn{3}{c|}{Triaxial} & \multicolumn{3}{c|}{Axisymmetric} \\ \hline 
\leavespace 
$\Delta t$ & $F^{\rm Eul}_{\rm def}$ & $F^{\rm Eul}_{\rm orth}$ & $F^{\rm Eul}_{\rm cor}$ 
	& $F^{\rm Eul}_{\rm def}$ & $F^{\rm Eul}_{\rm orth}$ & $F^{\rm Eul}_{\rm cor}$ \\ \hline 
\leavespace
$10$ & $9.04 \ 10^{-2}$ & $1.81 \ 10^{-2}$ & $4.56 \ 10^{-3}$
        & $7.66 \ 10^{-3}$ & $8.37 \ 10^{-4}$ & $0.00 $  \\
\leavespace
$1$ & $8.52 \ 10^{-2}$ & $8.54 \ 10^{-3}$ & $3.68 \ 10^{-5}$
        & $1.79 \ 10^{-3}$ & $1.13 \ 10^{-4}$ & $0.00$  \\
\leavespace
\leavespace
$.1$ & $9.06 \ 10^{-3}$ & $1.44 \ 10^{-3}$ & $3.55 \ 10^{-7}$
        & $2.51 \ 10^{-4}$ & $1.17 \ 10^{-5}$ & $0.00$  \\
\leavespace
$.01$ & $7.09 \ 10^{-4}$ & $1.62 \ 10^{-4}$ & $3.57 \ 10^{-9}$
        & $2.96 \ 10^{-5}$ & $1.18 \ 10^{-6}$ & $0.00$  \\
\hline
\end{tabular}
\end{center}
\caption{Average global energy errors over ten sample runs with randomly 
generated initial conditions and inertia tensors, integrated over the interval 
[0, 100] using versions of the forward Euler method.\label{random_errors}}
\end{table}
In table \ref{random_errors} we provide the average maximum errors in the 
energy for time steps $\Delta t= 10$, $1$, $.1$, and $.01$, using
for ten randomly generated initial conditions and inertia tensors each for
triaxial and axisymmetric bodies.

The separatrix is exactly captured if the infinitesimal updates 
$F^{\rm Eul}_{\rm orth}$ or $F^{\rm Eul}_{\rm cor}$, which coincide on the separatrix, are used. On 
the other hand, when $F^{\rm Eul}_{\rm def}$ was used to integrate ten sample trajectories 
with initial conditions at random points on the separatrices of rigid 
bodies with randomly generated inertia tensors, the average errors over
the integration interval [0, 500] were: $9.72 \ 10^{-2}$ for $\Delta t = 1$,
$3.89 \ 10^{-2}$ for $\Delta t = .1$, and $1.94 \ 10^{-3}$ for $\Delta t = .01$.

In the axisymmetric case, the forward Euler method with infinitesimal update
$F^{\rm Eul}_{\rm cor}$ associated to second order orbit approximation yields
the exact solution when the true exponential map is used as the algorithmic 
exponential. (If
the Cayley transform is used as the algorithmic exponential, then the 
orbits are captured exactly, but the algorithmic trajectories differ from 
the true trajectories by a time reparametrization.) Note that the `default'
generator and the orthogonal generator yield only first order orbit 
approximations even in the axisymmetric case. 

As implemented in our {\sl Mathematica} code, the version of the forward
Euler method with orthogonal 
algorithmic velocity is approximately 10\% slower than the naive version,
while the version that captures orbits to second order is approximately 
30\% slower than the naive version.

  \subsection{ Higher order methods for rigid body dynamics}
\begin{table}[h]
\begin{center}
\begin{tabular}{| r | c | c | c | c |} \hline
\leavespace
& \multicolumn{4}{c|}{Triaxial} \\ \hline 
\leavespace 
$\Delta t$ & $F^{RK2}_{\rm def}$ & $F^{RK2}_{\rm orth}$ & $F^{RK2}_{\rm dcor}$ & $F^{RK2}_{\rm ocor}$
\\ \hline 
\leavespace
$10$ & $7.22 \ 10^{-2}$ & $1.22 \ 10^{-2}$ & $6.11 \ 10^{-2}$ & $4.99 \ 10^{-2}$ \\
\leavespace
$1$ & $4.36 \ 10^{-3}$ & $1.06 \ 10^{-4}$ & $9.61 \ 10^{-4}$ & $5.67 \ 10^{-7}$ \\
\leavespace
$.1$ & $5.03 \ 10^{-6}$ & $1.10 \ 10^{-7}$ & $1.11 \ 10^{-8}$ & $2.14 \ 10^{-11}$\\
\leavespace
$.01$ & $5.01 \ 10^{-9}$ & $1.10 \ 10^{-10}$ & $2.94 \ 10^{-13}$ & $9.07 \ 10^{-15}$ \\ 
\hline
\leavespace
& \multicolumn{4}{c|}{Axisymmetric} \\ \hline
\leavespace
$\Delta t$ & $F^{RK2}_{\rm def}$ & $F^{RK2}_{\rm orth}$ & $F^{RK2}_{\rm dcor}$ & $F^{RK2}_{\rm ocor}$
\\ \hline
\leavespace
$10$ & $6.57 \ 10^{-3}$ & $2.27 \ 10^{-5}$ & $5054 \ 10^{-3}$ & $8.37 \ 10^{-7}$ \\
\leavespace
$1$ & $1.22 \ 10^{-4}$ & $2.44 \ 10^{-8}$ & $1.22 \ 10^{-5}$ & $9.09 \ 10^{-12}$ \\
\leavespace
$.1$ & $1.75 \ 10^{-7}$ & $2.46 \ 10^{-11}$ & $1.27 \ 10^{-10}$ & $3.15 \ 10^{-16}$\\
\leavespace
$.01$ & $1.75 \ 10^{-10}$ & $2.74 \ 10^{-14}$ & $6.61 \ 10^{-16}$ & $3.79 \ 10^{-15}$ \\ 
\hline
\end{tabular}
\end{center}
\caption{Average global energy errors over ten sample runs with randomly 
generated initial conditions and inertia tensors, integrated over the interval 
[0, 100] using versions of the Heun method.\label{random_errors2}}
\end{table}
\begin{table}[h]
\begin{center}
\begin{tabular}{| r | c | c | c | c | c | c |} \hline
\leavespace
& \multicolumn{6}{c|}{Triaxial} \\ \hline
\leavespace
$\Delta t$ & $F^{\rm sy4}_{\rm def}$ & $F^{\rm sy4}_{\rm orth}$ & $F^{\rm sy4}_{\rm dcor}$
& $F^{\rm sy4}_{\rm ocor}$ & $F^{RK4}_{\rm def}$ & $F^{RK4}_{\rm orth}$ \\ \hline
\leavespace
$10$ & $5.72 \ 10^{-3}$ & $3.37 \ 10^{-2}$ & $6.06 \ 10^{-3}$ & $1.90 \ 10^{-2}$
& $6.07 \ 10^{-2}$ & $3.59 \ 10^{-3}$ \\
\leavespace
$1$ & $3.62 \ 10^{-6}$ & $2.72 \ 10^{-7}$ & $1.55 \ 10^{-6}$ & $1.68 \ 10^{-7}$
& $2.15 \ 10^{-5}$ & $1.46 \ 10^{-7}$ \\
\leavespace
$.1$ & $3.56 \ 10^{-10}$ & $7.45 \ 10^{-12}$ & $3.61 \ 10^{-10}$
& $1.70 \ 10^{-12}$ & $2.97 \ 10^{-10}$ & $3.06 \ 10^{-12}$ \\ \hline
\leavespace
& \multicolumn{6}{c|}{Axisymmetric} \\ \hline
\leavespace
$\Delta t$ & $F^{\rm sy4}_{\rm def}$ & $F^{\rm sy4}_{\rm orth}$ & $F^{\rm sy4}_{\rm dcor}$
& $F^{\rm sy4}_{\rm ocor}$ & $F^{RK4}_{\rm def}$ & $F^{RK4}_{\rm orth}$ \\ \hline
\leavespace
$10$ & $2.47 \ 10^{-4}$ & $1.24 \ 10^{-7}$ & $2.47 \ 10^{-4}$  & $1.24 \ 10^{-7}$
& $7.45 \ 10^{-4}$ & $4.41 \ 10^{-8}$ \\
\leavespace
$1$ & $2.06 \ 10^{-8}$ & $1.26 \ 10^{-12}$ & $2.06 \ 10^{-8}$  & $1.26 \ 10^{-12}$
& $3.20 \ 10^{-7}$ & $4.47 \ 10^{-13}$ \\
\leavespace
$.1$ & $2.12 \ 10^{-13}$ & $0.00$ & $2.12 \ 10^{-13}$ & $0.00$
& $3.47 \ 10^{-12}$ & $3.22 \ 10^{-16}$\\ \hline

\end{tabular}
\end{center}
\caption{Average global energy errors over ten sample runs with randomly 
generated initial conditions and inertia tensors, integrated over the interval 
[0, 100] using several fourth order methods.\label{random_errors3}}
\end{table}

We implemented the four different versions of the Heun method given in section 
\ref{second_ordersection} for the rigid body system. 
Although the Heun methods $F_{\rm def}^{RK2}$ and $F_{\rm orth}^{RK2}$
described here are only second order accurate, the (local) discretization error 
in the energy is fourth order in the time step. 

The isotropy corrected versions used here take the form
\[
F^{RK2}_{\rm dcor}(m, \Delta t) := F_{\rm def}^{RK2}(m, \Delta t) + \Delta t^3 \, \sigma_{\rm def}(m) m,
\]
where
\[
\sigma_{\rm def}(m) := \frac {\left \langle{J}^n, 
	{u}\right \rangle}
	{\left \langle{J}^d, {u} \right \rangle}, \
\qquad \mbox{with} \qquad
\left\{ \begin{array}{l}
u_j := (m_k m_\ell)^2 \smallskip \\
J^n_j := - I_j (I_k + I_\ell)(I_k - I_\ell)^2 \smallskip \\
J^d_j := 4 \, I_1 I_2 I_3 I_j^2 (I_k - I_\ell)^2 
\end{array} \right .
\]
for any cyclic permutation $(j, k, \ell)$ of $(1, 2, 3)$ and
\[
F^{RK2}_{\rm ocor}(m, \Delta t) := F_{\rm orth}^{RK2}(m, \Delta t) + \Delta t^3 \, \sigma_{\rm orth}(m) m,
\]
where $\sigma_{\rm orth}$ is another, significantly more complicated, rational 
function in $m$ and the components of the inertia tensor.
The isotropy corrections $\sigma_{\rm def}$ and $\sigma_{\rm orth}$ 
given above determine algorithms yielding fourth order energy capture. 
If the body is axisymmetric, $F^{RK2}_{\rm dcor}$ preserves the energy to
fifth order. As Table \ref{random_errors2} shows, some of these algorithms 
appear to have better {\em global} energy capture than the single step 
discretization energy error analysis (which we carried out symbolically 
using {\sl Mathematica}) would suggest. Plots of the energy errors in
the sample integrations, with randomly generated initial conditions and
inertia tensors, show that the energy oscillates about a very slow drift 
away from the correct value.  
Note that the energy correction term for the fourth order symbolic expansion
method using the orthogonal generator is identically zero if the body is 
axisymmetric; hence the results generated by
$F^{\rm sy4}_{\rm orth}$ and $F^{\rm sy4}_{\rm ocor}$ coincide in this case.

We consider six fourth order geometric methods. Four utilize a series expansion
for the generator along a solution curve, while the other two use the RKMK4
algorithm (with the Cayley transform as the algorithmic exponential). 
Using the Cayley transform, the map $F^{\rm sy4}_{\rm def}$ determined by the default
generator for the rigid body system on $S^2$ is given by
\begin{eqnarray*}
\lefteqn{F^{\rm sy4}_{\rm def}(m, \Delta t) = \sum_{j = 1}^4 \frac {\Delta t^j}{j!} A^{(j - 1)}(m) 
        + \frac {\Delta t^3} {12} \left(  \| A(m)\|^2 A(m) + \dot A(m) \times A(m) \right)} \\
        && \qquad + {} \frac {\Delta t^4} {4!} \left(   \| A(m)\|^2 \dot A(m) 
	+ \ddot A(m) \times A(m) 
        + 2 \left \langle\dot A(m), A(m)\right \rangle A(m) \right),
\end{eqnarray*} 
where $A^{(j)}(m) = \frac{\partial^j }{\partial t^j}A({\mathcal F}_t(m))|_{t = 0}$.
The corresponding algorithm for the rigid body using the orthogonal generator is
\begin{eqnarray*}
F^{\rm sy4}_{\rm orth}(m, \Delta t) 
&=& \sum_{j = 1}^4 \frac {\Delta t^j}{j!} A_o^{(j - 1)}(m) \\
&& \quad + {} \frac {\Delta t^3} {12} \left(  \| A_o(m)\|^2 A_o(m) 
+ \left \langle\ddot A_o(m), m\right \rangle m \right) \\
&& \quad + {} \frac {\Delta t^4} {8} 
\left \langle\dot A_o(m), A_o(m)\right \rangle A_o(m).
\end{eqnarray*} 
The infinitesimal updates $F^{\rm sy4}_{\rm def}$ and $F^{\rm sy4}_{\rm orth}$ 
can be modified by the addition of an appropriate multiple of the argument $m$
to yield an additional order of energy, and hence orbit, capture. The
scalar correction functions, which are rational functions of the components of
$m$ and the inertia tensor, were determined by symbolic calculation.

  \section{ Application of geometric integration to numerical micromagnetics}

\label{num_results}

In this section, the geometric integrators developed in section
  \ref{sectionalgo} are applied to the Landau-Lifshitz-Gilbert
  equations of micromagnetics. The exact solution of this system is
  typically not available to us; indeed, it is the lack of precise
  analytical results for comparison which makes numerical
  micromagnetics a challenging field. In our examples, we chose the
  largest possible time steps for a given method that would lead the
  system to the solution computed by a higher order method (within prescribed tolerance).
  
As mentioned earlier, numerical micromagnetics has attracted much attention in the mathematical community, for several reasons. In this paper we are focussing on the time-stepping aspect of the problem. The application of geometric integration techniques in this context is relatively new, see for example, \cite{nigamlewis2000,krishna}. 
Recently, another technique which modifies existing integrators was developed for numerical micromagnetics, \cite{weinan}. This new integrator is of the ``step-and-project'' class, but is stable. 

  \begin{figure}[ht]
\centerline{\includegraphics[width=2.0in,height=2.5in]{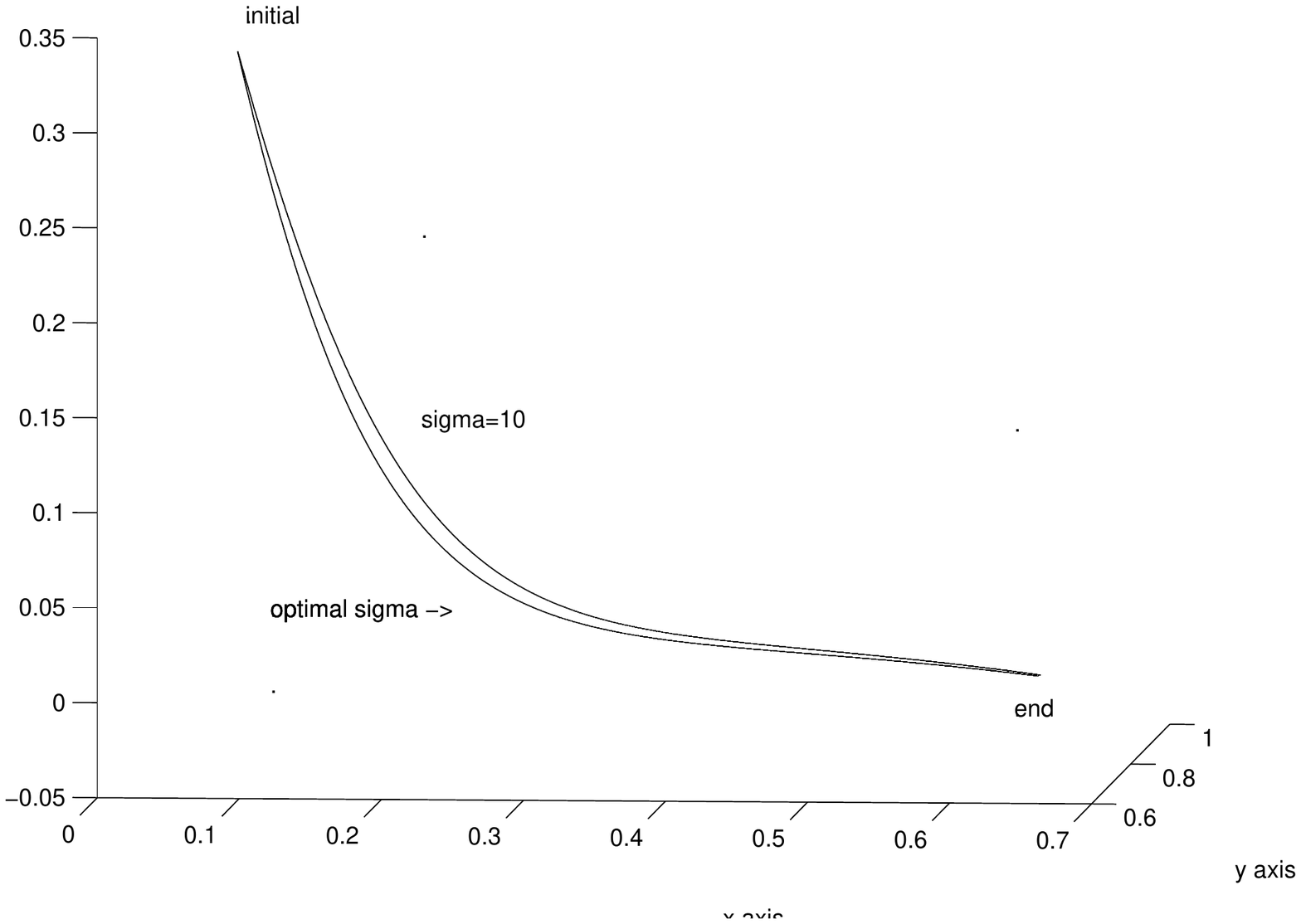} \hspace{0.4in}
\includegraphics[width=2.0in,height=2.5in]{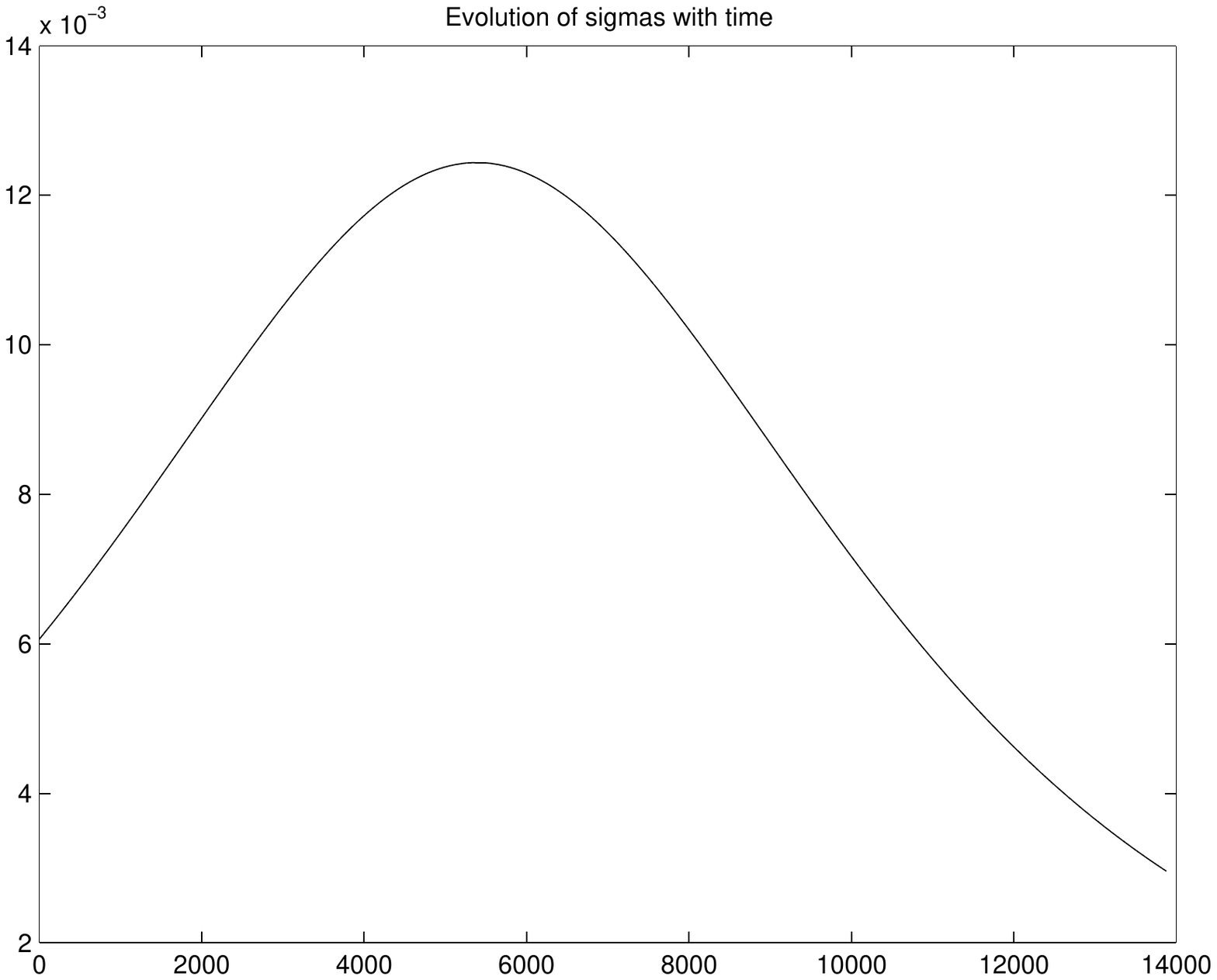}}                                                      
\caption{\small \it
The trajectories followed by the usual forward Euler $\Delta t=0.0001$
  and with the geometric forward Euler ($\Delta t=0.01$) with optimal
  $\sigma$ are almost identical; if we assume $\sigma=10$ (an arbitrary choice), the trajectory
  precesses more before reaching the final point. The figure on the
  right shows the evolution of the optimal $\sigma$. The applied field
  is uniform and weak,  specifically, $\bH_{\rm app} =(0.05,0.05, 0)$.}  \label{fig:1}  
\end{figure}
  \begin{figure}[ht]
\vspace{30mm}
\centerline{\includegraphics[width=2in,height=2.5in]{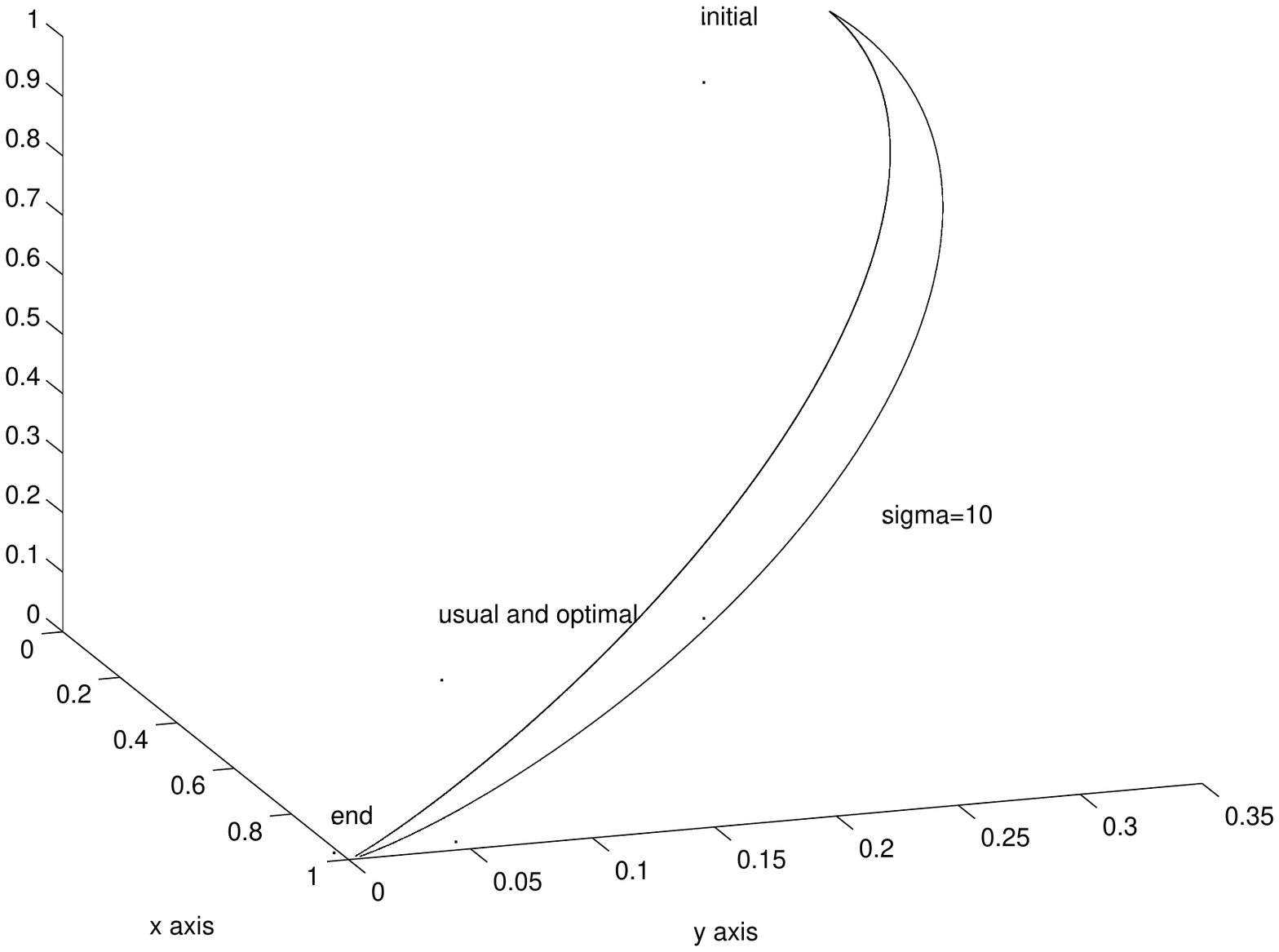}\hglue1in
\includegraphics[width=2in,height=2.5in]{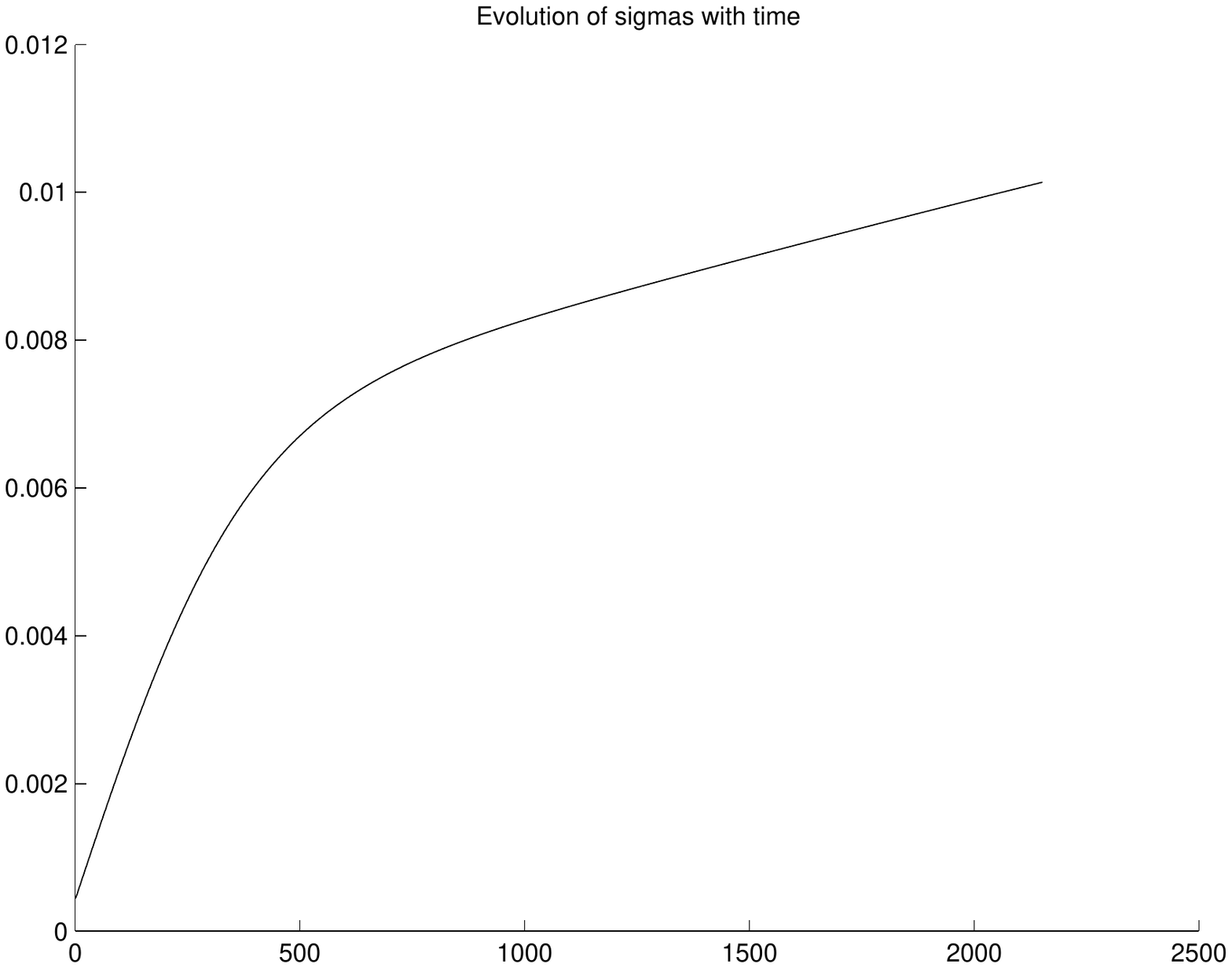} }                               
\caption{\small\it The evolution of one point in the ferromagnetic sample.
The trajectories followed by the usual forward Euler method with time step $\Delta t=0.0001$
  and by the geometric forward Euler method with optimal $\sigma$ and time step
  $\Delta t=0.01$ are almost identical; if we set $\sigma \equiv 10$ (an arbitrary choice) in the
  geometric forward Euler method, the trajectory
  precesses more before reaching the final point.  The figure on the
  right shows the evolution of the optimal $\sigma$. The applied field
  is uniform, $\bH_{\rm app}=(5,0, 0)$.} \label{fig:2}
\end{figure}
  
  We were particularly interested in the behavior of the free
  parameter $\sigma$ which appears in the geometric time-stepping
  algorithm (\ref{geom_eqs}). In the rigid body case  the
  parameter can be chosen to improve energy conservation. Here the
  system is dissipative and a criterion for the selection of $\sigma$ is
  not immediately obvious. For the forward Euler implementation, we
  can derive a relatively simple expression for a function $\sigma$ 
  that minimizes the discretization error. For higher order methods,
  analogous functions can be described in terms of the series expansions
  of the true and algorithmic flows, but the cost of computing these
  expansions, particularly for systems such as the LLG equations, rapidly
  becomes prohibitive. Work is in progress to determine computationally
  tractable criteria for the selection of the isotropy component for
  higher order methods. 

\subsection{Description of the model problem}
We describe a model problem for the LLG, for which the analytical solution was particularly simple. Recall that the LLG  for the magnetization $\bmu(x,t)$  is given by (\ref{LLG}), which we recall here for convenience:
\begin{equation}  \frac{\partial}{\partial t}\bmu=-\bmu\times \bH_{\rm eff}(\bmu)-\lambda \bmu\times(\bmu\times \bH_{\rm eff}(\bmu)), \qquad \|\bmu({\mathbf{x}})\|=1\,\,\,\forall \ {\mathbf{x}}\in \mathcal{B}.\end{equation}
 \begin{equation}\bH_{\rm eff}(\bmu)=A\Delta \bmu+\mu_0
  \left(-\nabla \phi+ \bH_{\rm app} \right)+K(\bmu\cdot \mathbf{e})\mathbf{e}.
  \end{equation} 

In the experiments that follow, we set $\mu_0=K=A=1.0$ and vary the applied field. These parameter values are not taken from actual physical data, and were selected solely for purposes of illustration. The saturation magnetization was $\|\bmu(\bx)\|\equiv 1$. 

We wish to construct a one-dimensional example where the computation of the demagnetizing energy would be simple. To this end, we assume the sample is contained in the infinite slab $\{(x,y,z)\,\vert 0\leq x \leq 1, y, z,\in \Rset^1\} $. We assume the magnetization $\bmu=\bmu(x)$, ie., the only variation in the magnetization is along the x-direction. Therefore, $\nabla \cdot \bmu = (\frac{\partial}{\partial x} \mu_1, 0, 0).$ We assume that there are 100 individual spins uniformly distributed along $x\,\in\,[0,1]$. These spins interact with each other through the exchange and demagnetizing fields. 

This example is admittedly a simplistic one; the true equilibrium solution for it can easily be found using analytical techniques. Therefore, the stopping criterion used was a comparison with the exact final equilibrium point. We see that the geometric integrators take trajectories which respect the point-wise constraints on the magnetization; conventional integrators do not. Thus, the paths traversed by these integrators will be different, as is seen in figures \ref{fig:1} and \ref{fig:2}. As the step-size is shrunk more and more, the trajectories will converge.

  \subsection{A first order method for the LLG}
 In the first set of numerical experiments, we implemented the
  geometrical analog of the forward Euler algorithm for the LLG equation.  We then tracked the 
evolution of the parameter $\sigma$ given by (\ref{movingsigma}), an
  expression derived through arguments of discretization error minimization.
We approximated the acceleration $\ddot \bM_n$ of trajectories $\bM(t)$ of (\ref{discretesystem}) using a 
  one-sided discrete approximation of the derivative of 
$\dot \bM_n=\mathbf{A}(\bM_n)\times \bM_n$.

   Figures \ref{fig:1} and \ref{fig:2} describe the trajectories followed by the
  over-damped LLG system (without the Larmor precession term) for two different
  applied fields.To address the issue of overall computational expense, we ran both a geometric forward Euler algorithm and the conventional forward Euler algorithm on the example introduced above. For each algorithm, we decreased the time step $\Delta t$ until the trajectories converged to within $5\%$.  We also required that the final equilibrium point corresponded to that computed by a fourth-order Runge-Kutta
  method with time step $\Delta t=0.0001$, to within a relative error of $1\%$.
The geometric forward Euler method yielded trajectories which converged, for this example, with time steps of  $\Delta t=0.01$ and a CPU 
time of 1.77 seconds. The usual forward Euler required a time step $\Delta t=0.0001$, with a CPU
time of 3.88 seconds,  to get similar behaviour.  In addition, while the pointwise norm of $\bM$ was conserved to machine accuracy by the geometric integrator, the standard forward Euler algorithm caused $\|\bM\|$ to increase to 1.001183806 times its usual value by the end of the run.  As a consequence, the trajectories traversed by the geometric and the usual algorithms differed, though they ended at the nearly same place. As we are only interested in the
  final equilibrium state of the system, we see the obvious merit of
  using the geometric integrator --- we can obtain accurate final states while
  using much larger time steps.
  
  We see the effect of varying the scalar
  functions $\sigma$ on the trajectories is that of changing the amount of precession in the trajectory.  We notice certain trends in the optimally chosen  function $\sigma(t)$ in
  figures \ref{fig:1} and  \ref{fig:2}, and  we shall
  investigate the relationship of these trends to the physical
  processes occurring at the same times in future work. 

  \begin{figure}[ht!]
\centerline{\includegraphics[width=2in,height=2.5in]{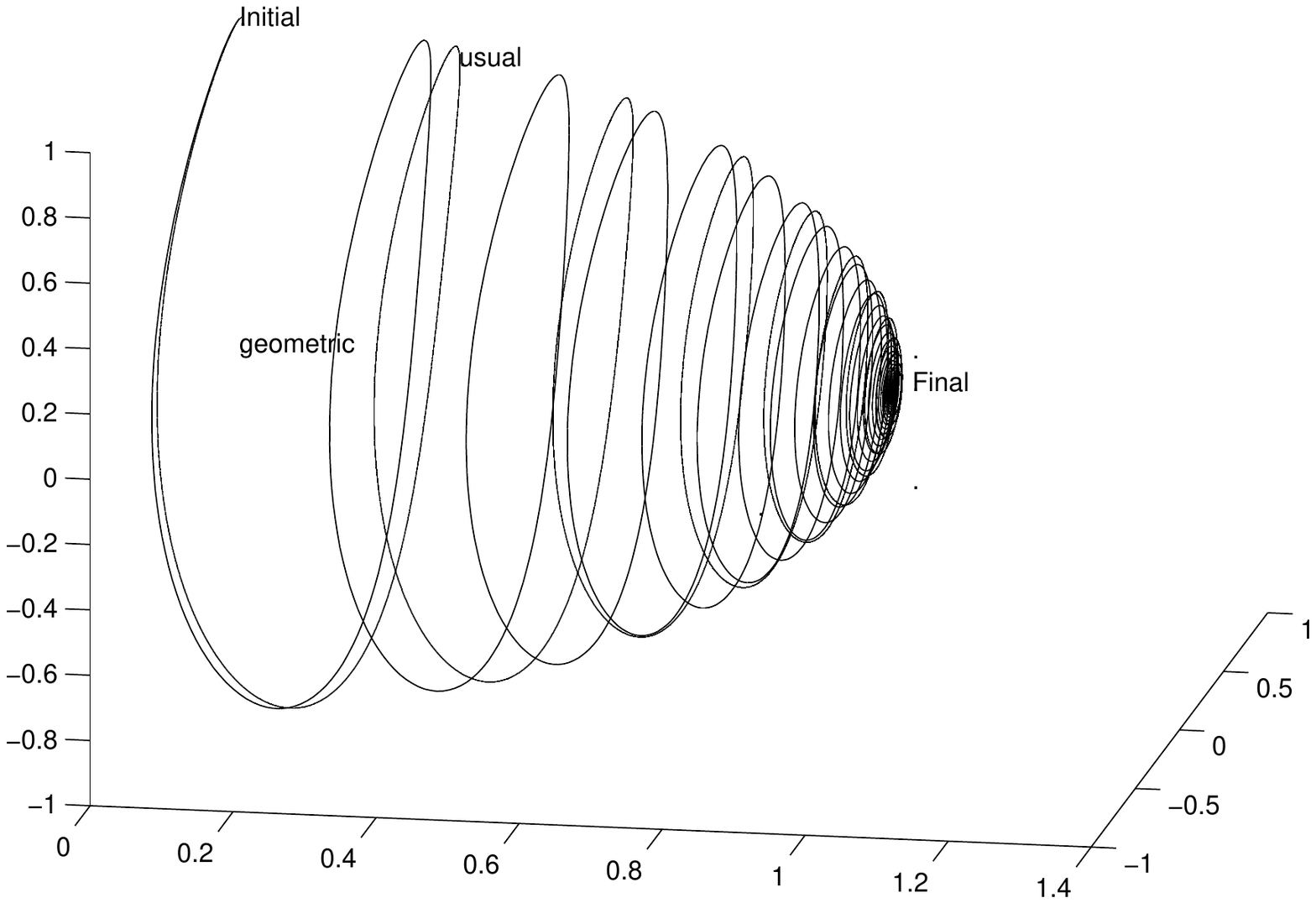}
\hglue1in\includegraphics[width=2in,height=2.5in]{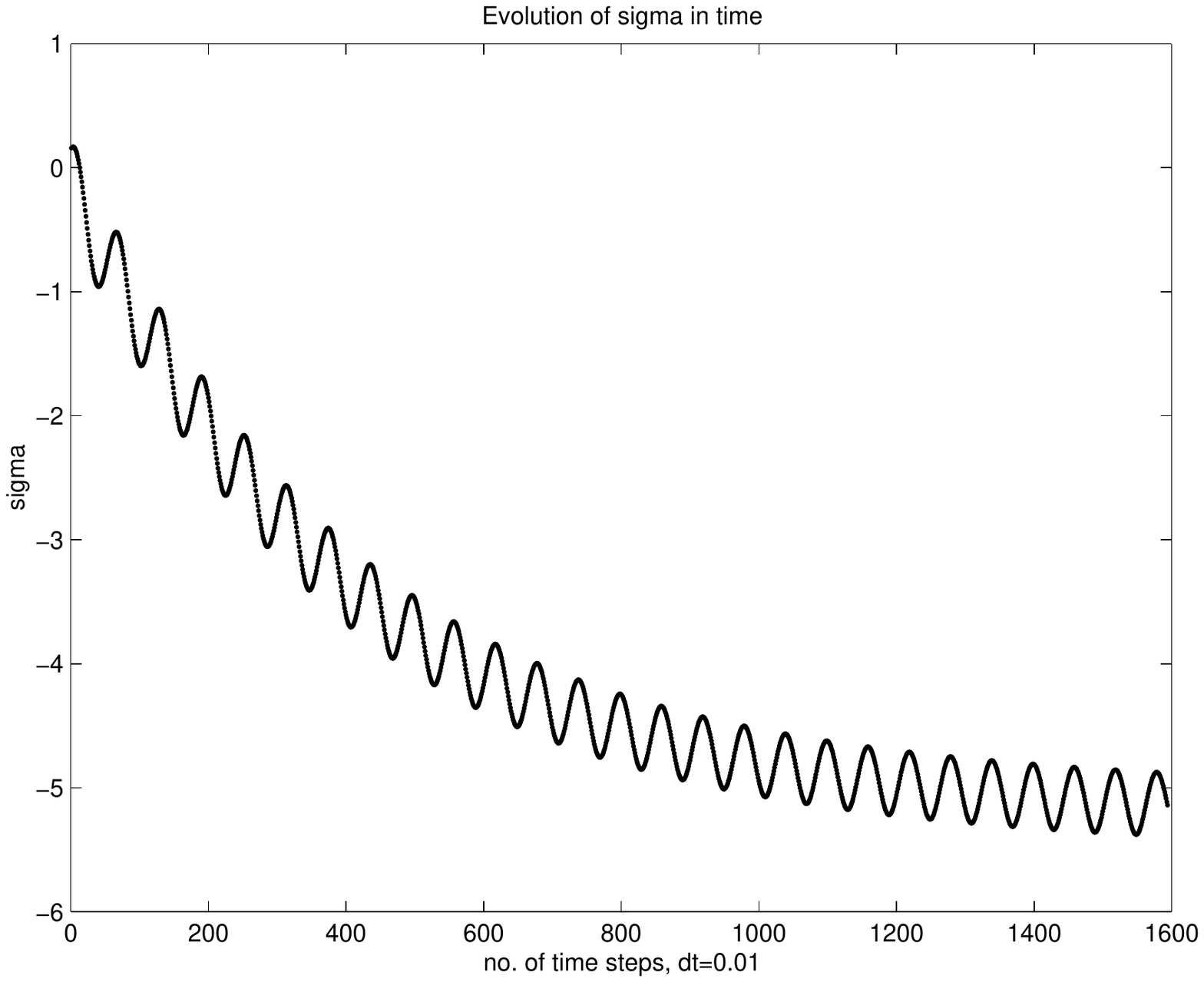}}
 \caption{\it \small  An example of the full LLG system, 
 with uniform applied field $\bH_{\rm app} =(5,0, 0)$ and damping parameter
  $\lambda = 0.05$. The trajectories were computed using the usual forward 
Euler method with time step $\Delta t=0.0001$ and the geometric 
  forward Euler method with time step $\Delta t=0.01$. Here the
  trajectory taken by the geometric method differs appreciably from
  those of the usual method, though the final states appear to be
  similar. The lefthand plot shows the evolution of one point in the
ferromagnetic sample; the righthand plot shows the evolution of the 
optimal sigma at that point.} \label{fig:4} 
  \end{figure}

\begin{figure}[ht!]
\vspace{30mm}
\centerline{\includegraphics[width=3in,height=3in]{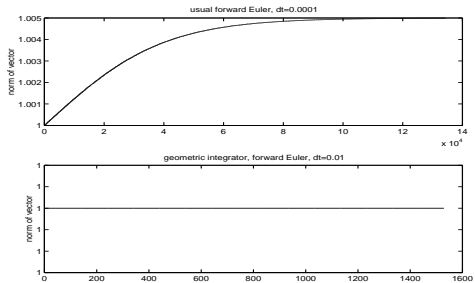}}
  \caption{\small \it Norm of magnetization in Fig(3). The geometric integrator exactly preserves
  the norm, even with a time step of $0.01$. The usual forward
  Euler method shows a drift in norm, even with a
  time step of $0.0001$. }   \label{fig:5}
\end{figure}

  In figure \ref{fig:4}, we implemented the code for the full LLG system,
  including the Larmor precession. The applied field
  is uniform, $\bH_{\rm app} =(5,0, 0)$. The damping parameter
  $\lambda$ was set to a low value, specifically $\lambda = 0.05$. The trajectories followed by the usual forward Euler and  
  $\Delta t=0.0001$
  and by the geometric forward Euler  with $\Delta t=0.01$ and  optimal
  $\sigma$ diverge appreciably, yet end at the same final state. The drift of the norm is now clearly
  visible (see figure (\ref{fig:5})). The usual forward Euler trajectory moves off the unit
  sphere in the standard Euler integration, while the geometrically integrated one does not. We see that the
  optimal $\sigma$ now varies more (figure \ref{fig:4}b).  

  \subsection{ A second order method for micromagnetics}
 In the next set of experiments, we implemented the geometrical
  versions of the Heun algorithm, derived in section \ref{sectionalgo}. We did not
  have an analytical expression for the optimal choice of
  $\sigma$. Therefore, we ran the experiments for several constant
  values of this parameter, and computed the order of convergence of
  the algorithm in $\Delta t$. 
  
  The results were interesting, and rather striking. As 
  $\sigma$ is varied, the order of convergence changes for the naive choice of
  generator. 
 What should be noted is that the geometric algorithm
 appears to converge more rapidly than a conventional Heun method; the order of convergence  was   $O(\Delta t^{2+\delta})$, as was borne out in repeated
  experiments.
  
  The norm of $\bM$ is conserved to machine precision for both generators.
\begin{figure}[ht!]
\vspace{30mm}
\centerline{\includegraphics[width=3in,height=3in]{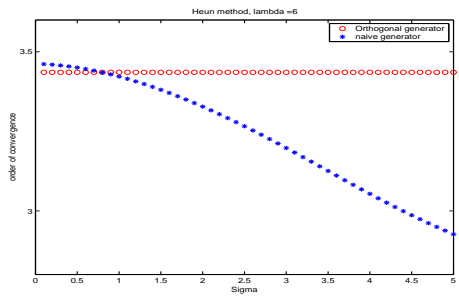}}
  \caption{\small \it Convergence orders of the Heun method for
  varying $\sigma$. We show experiments corresponding to two different
   generators.}   \label{fig:7}
\end{figure}

\subsection{A fourth order method for micromagnetics: RKMK4}
We now present experiments with a fourth order method derived in
  section \ref{sectionalgo}.  Lacking an analytical expression for the optimal choice of
  $\sigma$, we ran the experiments for varying constant
  values of this parameter, and computed the order of convergence of
  the algorithm in $\Delta t$. 
  
   \begin{figure}[ht!]
\vspace{30mm}
\centerline{\includegraphics[width=3in,height=3in]{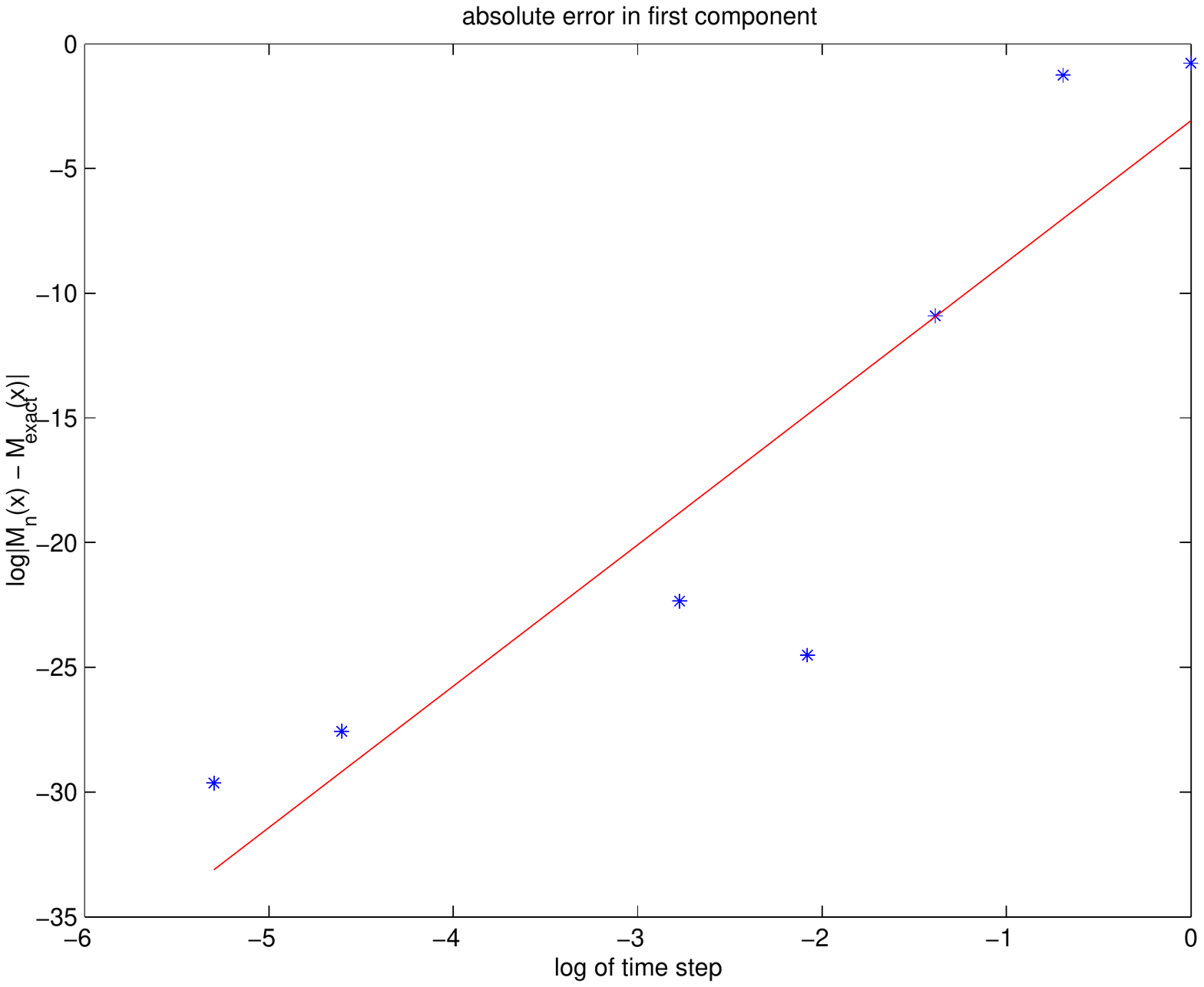}\includegraphics[width=3in,height=3in]{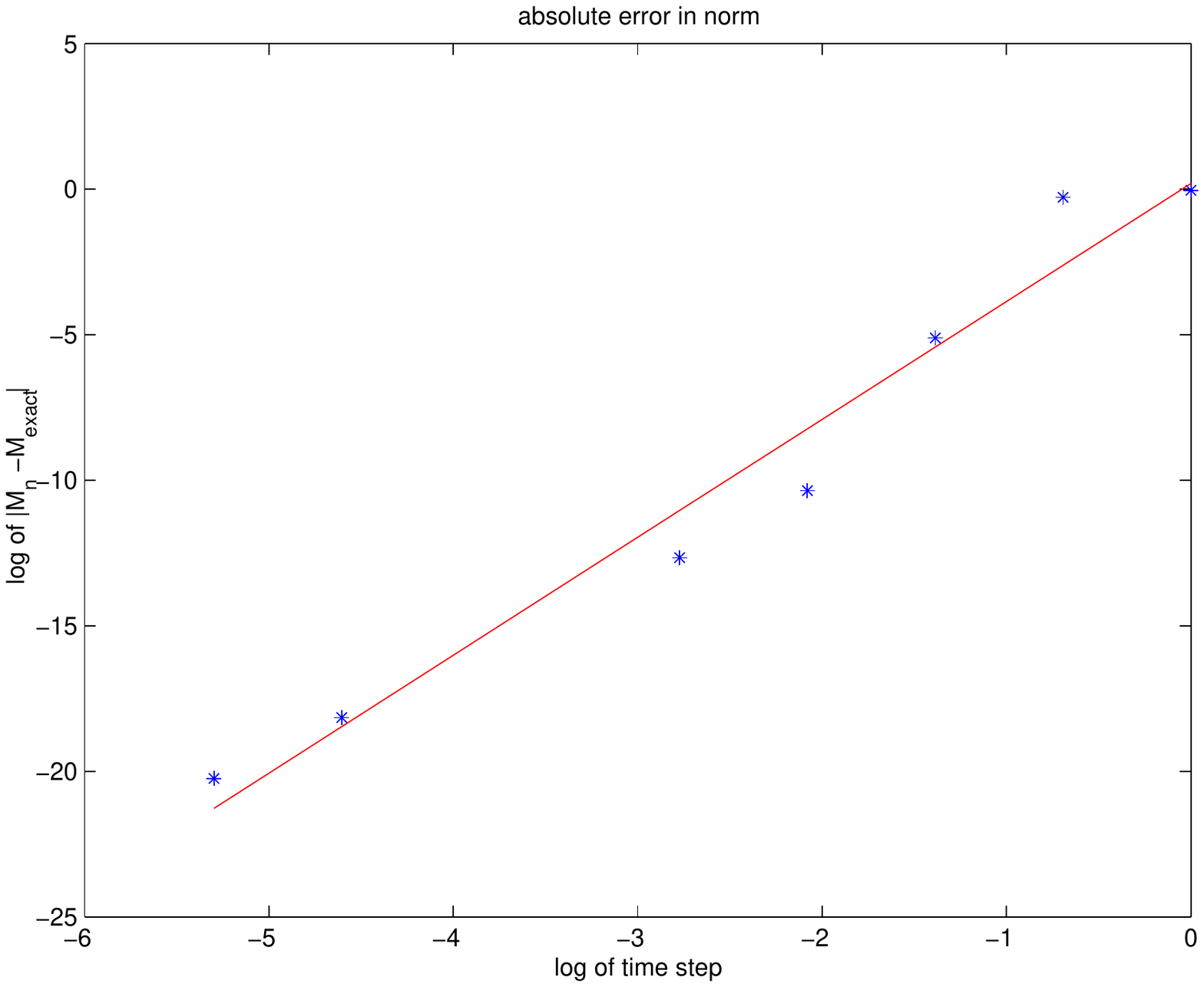}}
  \caption{\small \it On the left: Error in the first component of $\bM$ as a function
  of time step, $\sigma =0$. Here we see $O(\Delta t^{5.6})$
  convergence. On the right: Norm of the error of $\bM$ as a function
  of time step, with $\sigma=0$. Here we see $O(\Delta t^{4.03})$ convergence. }   \label{fig:7b}
\end{figure}
  
\begin{figure}[ht!]
\centerline{\includegraphics[width=3in,height=3in]{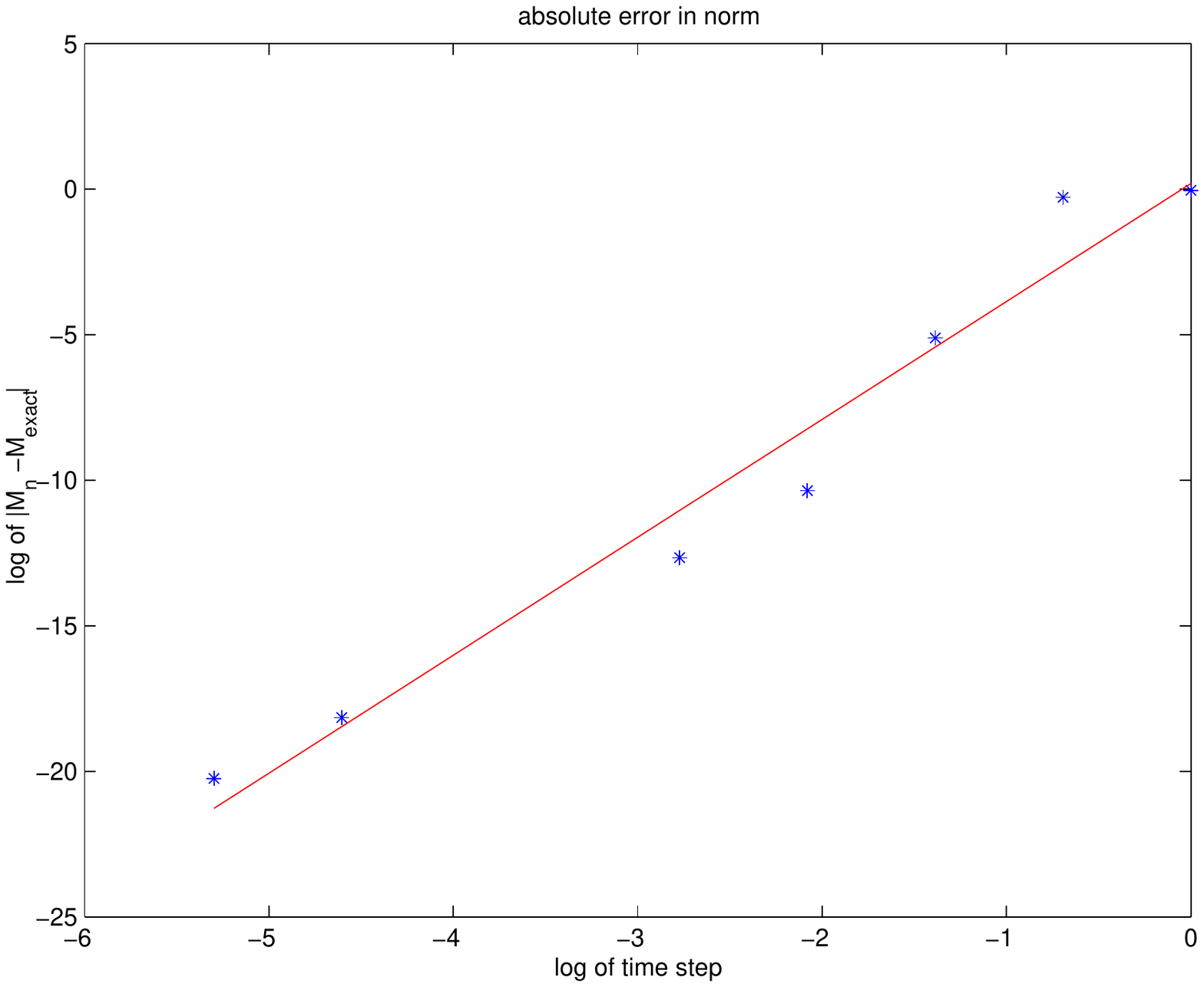}}
  \caption{\small \it Log(norm) of the magnetization over the
  integration interval }   \label{fig:9}
\end{figure}

In figure
  \ref{fig:9} we track $\|\bM\|$ over $[0,1]$ with a time step of $0.01$. The
  classical RK4 method without projection shows a drift in the norm;
  this drift is of the order of  $10^{-6}$, i.e., $O(\Delta t^3)$.  The Lie group
  integrator, on the other hand, shows no drift (up to machine
  precision).
  
\begin{figure}[ht!]
\centerline{\includegraphics[width=3in,height=3in]{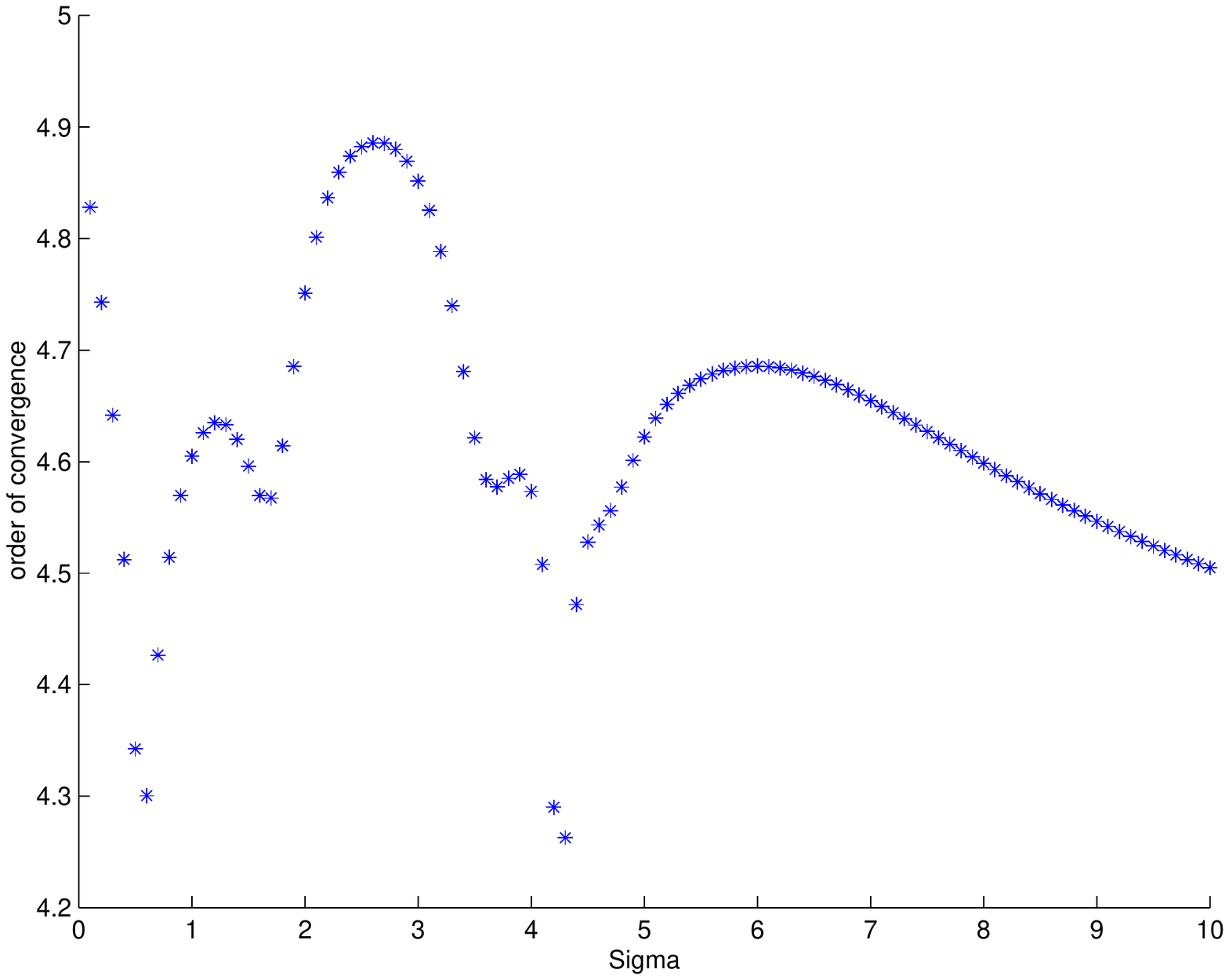}}
 \caption{\small \it Order of convergence of the algorithm as a
  function of $\sigma$, with damping parameter $\lambda =10.$}   \label{fig:10}
\end{figure}

 As we vary $\sigma$, we observe that the rate of
  convergence of the algorithm varies, see (figure \ref{fig:10}). Again, there is clearly some optimal
  value of this parameter. This behavior is even more pronounced for the RKMK4 
  method than for the Heun method.

{\bf Acknowledgements:} The authors would like to thank the referees for
their insightful and helpful suggestions. We also thank the Institute for Mathematics and its Applications, where the bulk of this research was conducted.

\end{document}